%% file: artrealhyp2.tex
\documentclass[a4paper,10pt]{article}
\usepackage[T1]{fontenc}
\usepackage{amsmath,amssymb,amscd,amsthm}
\usepackage[english]{babel}
\usepackage[latin1]{inputenc}
\usepackage{latexsym}
\usepackage[dvips]{graphicx}
\usepackage{epsfig}
\usepackage{graphics}
\usepackage{colortbl}
\usepackage{color}
\usepackage{latexsym}
\usepackage{ae}
\usepackage[cm]{aeguill}
\usepackage{vmargin}
\setmarginsrb{2 cm}{2 cm}{2 cm}{2 cm}{0 cm}{0 cm}{0,5 cm}{1 cm}
\begin{document}

\newtheorem{definition}{Definition}
\newtheorem{lemme}{Lemma}
\newtheorem{proposition}{Proposition}
\newtheorem{theoreme}{Theorem}
\newtheorem{conjecture}{Conjecture}
\newtheorem{corollaire}{Corollary}

\newcommand{\mb}{\mathbb}
\newcommand{\EM}{\ensuremath}
\newcommand{\norm}[1]{\EM{\left\| #1 \right\|}}

\title{Polyhedral realisation of hyperbolic metrics with conical singularities on compact surfaces}
\author{Fran\c{c}ois Fillastre}
\date{\today}
\maketitle

\selectlanguage{english}

\begin{abstract}

A Fuchsian polyhedron in hyperbolic space is a polyhedral surface
invariant under the action of a Fuchsian group of isometries (i.e. a group of isometries leaving globally invariant a
totally geodesic surface, on which it acts cocompactly). The induced metric on  a convex Fuchsian polyhedron is 
isometric to a hyperbolic metric with conical singularities of positive singular curvature on a compact surface of genus greater than one. We prove that these metrics are actually realised by exactly one convex Fuchsian polyhedron (up to global isometries). This extends a famous theorem of A.D. Alexandrov.
\end{abstract}
\tableofcontents

\section{Definitions and statements}

\paragraph{Metrics with conical singularities and convex polyhedra.}

Let $M_K^+$ be the simply connected (Riemannian) space of dimension $3$ of constant
curvature $K$, $K\in\{-1,0,1\}$. A \emph{convex polyhedron} is an intersection of
half-spaces of $M_K^+$. The number of half-spaces may be infinite, but the intersection is asked to be  locally finite:  each face must be a polygon with a finite number of vertices, and  the number of edges at each vertex must be finite.
A \emph{polyhedron} is a connected union of convex polyhedra. A \emph{polyhedral surface} is the boundary of a polyhedron and a \emph{convex polyhedral surface} is the boundary of a convex polyhedron.
A \emph{convex (polyhedral) cone}
in $M_K^+$ is a  convex polyhedral surface with only one vertex.
Note that the sum of the angles between the edges is strictly
between $0$ and $2\pi$.

A  \emph{metric of curvature $K$ with conical singularities with positive singular curvature} on a compact surface $S$ is a
(Riemannian) metric of constant curvature $K$ on $S$ minus $n$
points $(x_1,\ldots,x_n)$ such that the neighbourhood of each $x_i$
is isometric to the induced metric on the neighbourhood of the vertex of a convex
cone in $M_K^+$. The $x_i$ are called the \emph{singular points}. By
definition the set of singular points is discrete, hence finite since
the surface is compact.

An \emph{invariant polyhedral surface} is a pair $(P,F)$, where $P$ is a polyhedral surface in $M_K^+$ and 
$F$ a discrete group of isometries of $M_K^+$ such that $F(P)=P$ and $F$ acts freely on $P$. The group $F$ is called 
the \emph{acting group}.

If there exists an invariant polyhedral surface $(P,F)$ in $M_K^+$ such that the induced metric on $P/F$ is isometric to a metric $h$ of curvature $K$ with conical singularities   on a surface $S$, we say that $P$ \emph{realises} the metric $h$ (obviously the singular points of $h$ correspond to the vertices of
$P$, and $F$ is isomorphic to the fundamental group of $S$). In this case we say that $h$ is \emph{realised by a unique invariant polyhedral surface $(P,F)$} if $P$
is unique up to isometries of $M_K^+$.

Let $P$ be the boundary of a convex compact polyhedron in $M_K^+$. The
induced metric on $P$ is isometric to a  metric of constant curvature $K$ with conical singularities of positive singular curvature on
the sphere.

 A famous
theorem of A.D. Alexandrov asserts that each such  metric on the
sphere is realised by the boundary of a unique convex compact polyhedron of $M_K^+$ \cite{Aleks,Buse,Pogo} - in this case the acting group $F$ is the trivial one. 

In this paper we prove 
\begin{theoreme}\label{realisation}
A hyperbolic metric with conical singularities of positive singular curvature on a compact surface $S$ of genus $>1$ is realised by a unique convex Fuchsian polyhedron in hyperbolic space (up to global isometries).
\end{theoreme}

A \emph{Fuchsian polyhedron} is a   polyhedral surface  invariant under the action 
of a Fuchsian group of hyperbolic space $\mathbb{H}^3$. A \emph{Fuchsian group of hyperbolic space}  
is a discrete group of orientation-preserving isometries leaving globally invariant a
totally geodesic surface, on which it acts cocompactly and without fixed points. 
  The idea to use them comes from Gromov \cite{gromovdiffrel}. Analog statement can be found in \cite{Schpoly}, see further.

The general outline of the proof of Theorem \ref{realisation} is very classical and has been used in several other cases, starting from A.D. Alexandrov's works. Roughly speaking, the idea is to endow with   suitable topology both the space of convex Fuchsian polyhedra with $n$ vertices and the space of  corresponding metrics, and to show that the map from one to the other given by the induced metric is a homeomorphism.

The difficult step is to show the local injectivity of the map ``induced
metric''. This is equivalent to a statement on infinitesimal rigidity of convex
Fuchsian polyhedra. The  section \ref{rigidite}  is devoted to this result.

\paragraph{Example of convex Fuchsian polyhedra.}

Consider the
set of points at  constant (hyperbolic) distance of a totally geodesic surface
$P_{\mathbb{H}^2}$  and denote by  $M$ the subset of points which
are on one side of  $P_{\mathbb{H}^2}$. 

Obviously $M$ is
globally invariant under the action of any Fuchsian group $F$ 
leaving  $P_{\mathbb{H}^2}$ invariant. Moreover, $M$  has the properties of being strictly convex and umbilic (an \emph{umbilic surface} is a surface such that its principal curvatures are
the same for all points. There are other kinds of umbilic
surfaces in the hyperbolic space, but in all this text, the
expression ``umbilic surface'' means a complete surface at constant
distance of $P_{\mathbb{H}^2}$ and contained in one of the half-spaces 
bounded by $P_{\mathbb{H}^2}$).

Take $n$ points $(x_1,\ldots,x_n)$ on $M$, and let $F$ act on  these points. We denote by  $E$
the boundary  of the convex hull of the points  $f x_i$, for all $f\in F$ and $i=1\ldots n$. By construction, the convex polyhedral surface $E$ is globally invariant under the action $F$: it is a convex Fuchsian polyhedron.

%\paragraph{Immediate corollaries of Theorem \ref{realisation}.}
\paragraph{Global rigidity of convex Fuchsian polyhedron.} A
polyhedral surface is called \emph{globally rigid} if any polyhedral surface that is isometric to it is in fact congruent. A direct consequence of the
uniqueness of the convex Fuchsian polyhedron realising the induced
 metric is
\begin{theoreme}
Convex Fuchsian polyhedra in hyperbolic space are globally rigid among  convex Fuchsian polyhedra.
\end{theoreme}

\paragraph{Hyperbolic manifolds with polyhedral boundary.}

\begin{figure}[h]
\begin{center}\includegraphics{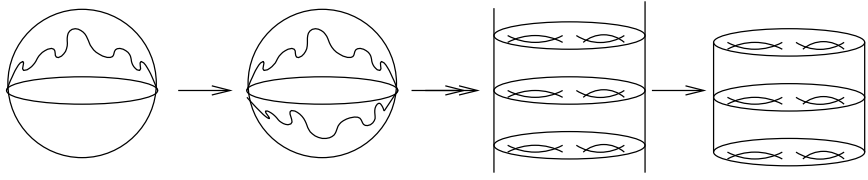}\end{center}
\caption{From a Fuchsian polyhedron to a Fuchsian manifold. \label{symetrie}}
\end{figure}

Take a convex Fuchsian polyhedron $(P,F)$ and consider the Fuchsian polyhedron
$(P',F)$ obtained by the reflection on the invariant surface $P_{\mathbb{H}^2}$ of $F$. Then cut the hyperbolic space
along $P$ and $P'$, and keep the  component bounded by $P$ and $P'$. The quotient of this manifold by $F$ is a kind of
hyperbolic manifold called \emph{Fuchsian manifold} (with
convex polyhedral boundary): these are compact hyperbolic manifolds with  boundary with an
isometric involution fixing a compact hyperbolic surface (the symmetry relative to
$P_{\mathbb{H}^2}/F$), see Figure \ref{symetrie}. In this case we obtain a Fuchsian manifold with convex polyhedral boundary, and all the Fuchsian manifolds with convex polyhedral boundary can be
obtained in this way: the lifting to the universal cover of  a component of the boundary of the Fuchsian manifold 
gives a convex Fuchsian polyhedron in the hyperbolic space. Then  Theorem \ref{realisation} says
exactly that for a choice of the metric $h$ on the boundary, there
exists a unique  metric on the manifold such that it is a Fuchsian manifold with convex polyhedral boundary and the
induced metric on the boundary is isometric to $h$:
\begin{theoreme}
The metric on a Fuchsian manifold with convex polyhedral boundary is determined by the induced metric on its boundary.
\end{theoreme}
This is a part of
\begin{conjecture}\label{conjecture variete}
Let $h$ be a hyperbolic metric on  a compact  manifold $M$  of dimension $3$  such that $\partial M$ is polyhedral and convex. Then the induced metric on $\partial M$ is a hyperbolic metric with conical singularities with positive singular curvature. Each hyperbolic metric with conical singularities with positive singular curvature on $\partial M$ is induced on $\partial M$ for a unique choice of $h$.
\end{conjecture}

The  statement of Conjecture \ref{conjecture variete} in the case where the boundary is smooth and strictly convex as been proved in \cite{Schconvex} (the existence part was found in \cite{lab2}). Remark that A.D. Alexandrov's theorem is a part of Conjecture \ref{conjecture variete} for
the case of the hyperbolic ball.

The smooth analog of Conjecture \ref{conjecture variete} provides a smooth version of Theorem \ref{realisation}. The existence part of this smooth statement was done in \cite{gromovdiffrel}.

\paragraph{Towards a general result.}

The Lorentzian space-forms of dimension $3$ are the de Sitter space $dS^3$ (with curvature $1$), the Minkowski space $\mathbb{R}^3_1$ (with curvature $0$) and the anti-de Sitter space $AdS^3$ (with curvature $-1$), see e.g. \cite{oneill}. We denote them by $M_K^-$, where $K$ is the curvature.

Such spaces contain surfaces on which the induced metric is Riemannian (think about the hyperbolic plane in the Minkowski space), and these surfaces are called \emph{space-like}. A space-like convex polyhedral cone in a Lorentzian space-form has a negative singular
curvature at its vertex (i.e. the sum of the angles between its edges is $>2\pi$).

A
theorem of Rivin--Hodgson \cite{rivinthese,RivHod} says that each metric of curvature $1$  on the sphere with conical singularities with negative singular curvature such that its closed geodesics have length $>2\pi$  is realised by a unique space-like convex   polyhedral surface homeomorphic to the sphere in the de Sitter space (beware that 
$\pi_2(dS^3)\not= 0$).

In \cite{Schpoly} it is stated that each metric of curvature $1$  on a compact surface of genus $>1$  with conical singularities with negative singular curvature such that its contractible geodesics have length $>2\pi$  is realised by a unique space-like convex Fuchsian    polyhedron  in the de Sitter space (the definition of a Fuchsian polyhedron in a Lorentzian space-form is the same as in the hyperbolic space, after replacing ``totally geodesic surface'' by ``umbilic hyperbolic surface''). 

We think that each constant curvature $K$ metric with conical singularities with constant sign singular curvature $\epsilon \in \{-,+\}$ on  a compact
surface can be realised in $M_K^{\epsilon}$ by a unique (space-like) convex  polyhedral surface invariant under the action of a representation of the fundamental group of the surface in a group of isometries of dimension $3$ - with a condition 
on contractible geodesics in the cases $K=1,\epsilon=-$. The images of these representations are trivial for genus $0$, parabolic for genus $1$ and Fuchsian for genus $>1$ (a parabolic isometry fixes a point on the boundary at infinity).

The present paper proves this assertion for  hyperbolic metrics with conical singularities with
positive singular curvature on compact surfaces of genus $>1$.

All the combinations with genus, curvature and sign of the singular curvature are not possible because of Gauss--Bonnet Formulas \cite{Troyanovarticle2}.
 The complete list of results would be, if $g$ is the genus of the compact surface (we imply that the polyhedral surfaces are (space-like) convex):
  \begin{itemize}
\item[$g=0$]
  \begin{itemize}
 \item[] $K=-1$, $\epsilon = +$: boundary of a compact polyhedron  in $\mathbb{H}^3$ (Alexandrov);
 \item[] $K=0$, $\epsilon = + $: boundary of a compact polyhedron in $\mathbb{R}^3$ (Alexandrov);
  \item[] $K=1$, 
  \begin{itemize} \item[] $\epsilon = + $: boundary of a compact polyhedron in $\mathbb{S}^3$ (Alexandrov);
   \item[] $\epsilon = - $ and length of the contractible  geodesics $>2\pi$: polyhedral surface homeomorphic to the sphere in $dS^3$ (Rivin--Hodgson);
 \end{itemize}
 \end{itemize}
\item[$g=1$] 
 \begin{itemize}
 \item[] $K=-1$, $\epsilon = +$: parabolic polyhedron in $\mathbb{H}^3$;
 \item[] $K=1$, $\epsilon = - $ and length of the contractible  geodesics $>2\pi$: parabolic polyhedron in $dS^3$;
 \end{itemize}
 \item[$g>1$]  
 \begin{itemize}
 \item[] $K= -1$,
 \begin{itemize}
 \item[] $\epsilon=+$: Fuchsian polyhedron in $\mathbb{H}^3$ (this paper);
 \item[] $\epsilon = -$: Fuchsian polyhedron in $AdS^3$;
  \end{itemize}
  \item[] $K=0$, $\epsilon =-$: Fuchsian  polyhedron in $\mathbb{R}^3_1$;
 \item[] $K=1$, $\epsilon = - $ and length of the contractible  geodesics $>2\pi$: Fuchsian polyhedron in $dS^3$ (Schlenker).
 \end{itemize}
 \end{itemize}
The proofs of the others cases for $g>1$ would be close to the one presented here \cite{artrealisationlorentz}.

\paragraph{Acknowledgements.}
The material in this paper is a part of my PhD thesis under the direction of
B. Colbois and J.-M. Schlenker. For that reason, they played a
crucial part in the working out of these results. I also want to
thank M. Troyanov for his useful comments.

\section{Infinitesimal rigidity}\label{rigidite}

\subsection{Background about infinitesimal isometric deformations}

A \emph{Killing field} of a constant curvature space $M_K^+$  is  a vector
field of $M_K^+$ such that the elements of its local 1-parameter group are
isometries (see e.g. \cite{GaHuLa}). An \emph{infinitesimal isometric deformation} of a
polyhedral surface consists of 
\begin{itemize}
\item a triangulation of the polyhedral surface given by a triangulation of each face, such that no new vertex arises,
\item a Killing field on each face of the
triangulation such that  two Killing fields on two adjacent triangles are equal on the common edge. 
\end{itemize}
The edges of the triangulation which were not edges of the polyhedral surface are called \emph{additional edges}.

%Then an infinitesimal isometric deformation of a polyhedral surface is uniquely
%determined by its values at the vertices of the polyhedral surface.

An infinitesimal isometric deformation is called \emph{trivial}
if it is the restriction to the polyhedral surface of a global Killing
field. If all the infinitesimal isometric deformations of a
polyhedral surface are trivial, then the polyhedral surface is said to be \emph{infinitesimally rigid}.

\subsubsection{Infinitesimal rigidity of polyhedral convex caps which may have infinitely vertices accumulating at the boundary}

A  \emph{polyhedral convex cap}  is a   convex polyhedral surface $C$, with boundary $\partial C$, in the Euclidean space homeomorphic to  the (closed) disc, such that $\partial C$ 
lies in a totally geodesic plane, and such that the orthogonal projection onto this plane is a bijection between $C$ and the domain of the plane inside $\partial C$ (up to global isometries, we suppose that $\partial C$ lies in the horizontal plane).

This paragraph is dedicated to the proof of:
\begin{proposition} \label{bordcon2}
  If the vertical component  of an infinitesimal isometric deformation of a polyhedral convex cap vanishes on the boundary, then the deformation is trivial.
\end{proposition}

Note that our definitions allow polyhedral convex caps with an infinite number of vertices which accumulate at the boundary, and also infinitesimal deformations which diverge on the boundary.

A smooth version of  Proposition \ref{bordcon2} has been known for a long time \cite[Thm 1 Chap IV § 7]{Pogo} - this reference contains also polyhedral results.

\begin{definition}
 Let $C$ be a polyhedral convex cap and $p$ a vertex of $C$. Consider the spherical polygon which is the intersection of 
 a triangulation of $C$ with a little sphere centred at $p$
(such that it intersects only edges incident with $p$). The \emph{link} of
$C$ at $p$ is the image of this polygon by a homothety sending the
little sphere to the unit sphere.
\end{definition}

As $C$ is convex, the link is a convex spherical polygon (with some interior angles may be equal to $\pi$ because of the additional edges of $C$).

We denote by $u$ the vertical component of an infinitesimal isometric
deformation of $C$. The definition of an infinitesimal isometric deformation implies that it is a continuous vector field defined 
on the interior of the convex cap. It follows that the function $u$ is  continuous 
 on the interior of the convex cap. 
 
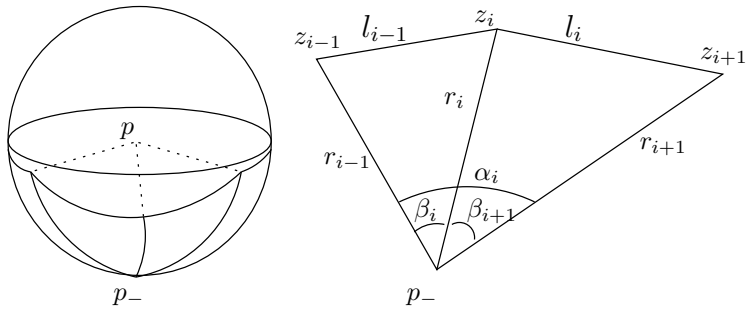
\begin{figure}\begin{center}
\input{angleconvexe.pstex_t}\end{center}
\caption{\label{conelink}Some notations for the link at $p$.}
\end{figure}

 Up to a Euclidean isometry, we can consider
that at a vertex $p$,  $u(p)=0$.  Then, if
the value of $u$ is positive (resp. negative) at a point on an edge joining $p$ to another vertex, it means that the vertical component of the
deformation at this point is greater (resp. less) than  at
$p$. This doesn't depend on the choice of the point on the edge. 
In particular, it is the case for the point which is sent by a 
homothety to a vertex of the link.
In this case, we say that the corresponding vertex of the link \emph{goes up} (resp.
\emph{goes down}).
\begin{proposition}\label{proposition1}
The link cannot go up or go down,  i.e. all the vertices of the
link cannot go up or go down at the same time.
\end{proposition}
We build a spherical polyhedral convex cone from the link by joining all the vertices of the link with the south pole - denoted by  $p_-$ - of the sphere containing  the link.
We denote by $(z_1,\ldots,z_n)$ the vertices of the link,
$\beta_i$ the angle at $p_{-}$ between the (spherical) segments $p_{-}z_i$ and $p_{-}z_{i+1}$,  $\alpha_i$
the angle at $p_{-}$ between the (spherical) segments $p_{-}z_i$ and $p_{-}z_{i+2}$ and $r_i$ the length of the (spherical) segment
between $p_{-}$ and $z_i$ (see Figure \ref{conelink}).

Without loss of generality, we consider that all the vertices of the link go up.
  By definition the infinitesimal isometric deformation is isometric on
the faces of the polyhedral convex cap, it implies that the lengths  $(l_{1},\ldots,l_{n})$ of the edges of the link don't
change under the deformation:  each $\beta_i$ is a function of
$r_{i-1}$ and $r_{i}$.

We recall the well known ``spherical law of cosines (for the sides)''.
Let $a,b,c$ be the length of the edges of a spherical triangle, and $\alpha$ the angle at the opposite vertex from the edge of length $a$. Then
\begin{equation}
\cos a=\cos b \cos c + \sin b \sin c \cos \alpha.\nonumber
\end{equation}

\begin{lemme}[{Corollary of the Cauchy Lemma, \cite[18.7.16]{berger5},\cite{sabitov}}]
If a (convex) quadrilateral of the sphere is deformed such that the lengths 
of the edges remain constant under the deformation and such that two opposite angles increase, then the two others angles decrease (an angle may be equal to $\pi$).
\end{lemme}
\begin{lemme}\label{fonction angle deformation}
If  $r_{i-1}$ and $r_{i+1}$ are fixed, then
\begin{equation}
\frac{\partial \beta_i}{\partial r_i}(r_{i-1},r_{i})+
\frac{\partial \beta_{i+1}}{\partial r_i}(r_{i},r_{i+1})<0.\nonumber
\end{equation}
\end{lemme}
\begin{proof}
We call $\theta$ the angle at $z_{i-1}$  
of the triangle
$(p_-,z_i,z_{i-1})$.  The  ``spherical law of cosines'' applied to this triangle says 
  that 
 \begin{equation}
 \cos r_i=\cos l_{i-1} \cos r_{i-1} +\sin l_{i-1} \sin r_{i-1} \cos \theta. \nonumber
 \end{equation}
  As $l_{i-1}$ and $r_{i-1}$ are supposed to be fixed, and as the sines are positive, we deduce from this formula that $\theta$  is a
strictly increasing  function of $r_{i}$.

In the same way, the angle at  $z_{i+1}$ of the triangle  $(p_-,z_{i+1}, z_i)$ is a strictly increasing  function of $r_{i}$.

In this case, the  corollary of the Cauchy Lemma  says that
the angles at $p_-$ and $z_i$ of the quadrilateral $(p_-, z_{i+1},
z_i, z_{i-1})$ decrease. The first of these angles was called
$\alpha_i$, and it is the sum of  $\beta_i$ and $\beta_{i+1}$. It
shows that this sum (strictly) decreases when $r_i$ (and only
$r_i$) increases, and this is another way to state the lemma.
\end{proof}
\begin{proof}[Proof of  Proposition \ref{proposition1}]
We consider $(\dot{r}_1,\ldots,\dot{r}_n)$, where $\dot{r}_{i}$ means
 $\frac{d}{dt}r_{i}(t)\vert_{t=0}$, a deformation (defined on each vertex) of the link such that
all the vertices go up, i.e. $\dot{r}_i>0$ $\forall i$.
As the sum of the angles $\beta_j$ goes around $p$ we have
\begin{equation}\sum_{j=1}^{n}\beta_j=2\pi\nonumber\end{equation}
and this remains true under the deformation:
 \begin{equation}\label{sommenulle}\sum_{j=1}^{n}\dot{\beta}_j=0.\end{equation}
But on other hand
\begin{equation}\dot{\beta}_j=\frac{\partial \beta_j}{\partial r_{j-1}}\dot{r}_{j-1}
+\frac{\partial \beta_{j}}{\partial
r_{j}}\dot{r}_{j},\nonumber\end{equation} and we get a contradiction from Lemma \ref{fonction angle deformation},  by
replacing in Equation (\ref{sommenulle}) (with a cyclic
notation $\beta_{n+1}=\beta_{1}$):
\begin{eqnarray}
\ 0&=&\sum_{j=1}^{n} \left( \frac{\partial \beta_{j}}{\partial r_{j-1}}\dot{r}_{j-1}+
\frac{\partial
  \beta_{j}}{\partial r_{j}}\dot{r}_{j}\right)  =\sum_{k=1}^{n}\frac{\partial
  \beta_{k}}{\partial r_{k-1}}\dot{r}_{k-1}+\sum_{i=1}^{n}\frac{\partial
  \beta_{i}}{\partial
  r_{i}}\dot{r}_{i} \nonumber \\
\ &\stackrel{i=k-1}{=}&\sum_{i=-1}^{n-1}\frac{\partial \beta_{i+1}}{\partial
  r_i}\dot{r}_i+\sum_{i=1}^{n}\frac{\partial \beta_{i}}{\partial
  r_{i}}\dot{r}_{i}
=\sum_{i=1}^{n}\left(\underbrace{\frac{\partial \beta_i}{\partial
  r_i}+\frac{\partial \beta_{i+1}}{\partial r_{i}}}_{<0}\right)\underbrace{\dot{r}_{i}}_{>0}.\nonumber
\end{eqnarray}
\end{proof}

\begin{corollaire}
The function $u$ does not attain a local extremum in the interior of the convex cap.
\end{corollaire}
\begin{proof}
  Suppose that a local extremum is reached at a vertex $p$ (as the deformation is isometric, it is clear that a local extremum can't be reached at a point on a face or an edge). That means that the values of  $u$ at points on the edges from $p$  are all greater
(or less) than $u(p)$. It implies that the link at  $p$ goes up (or down), that's impossible by the preceding proposition. 
\end{proof}
If we make the hypothesis that $u$ vanishes at the boundary, this
corollary implies that $u$ vanishes for all the vertices of the
polyhedral convex cap, and that proves the triviality of the
infinitesimal isometric deformation associated to $u$. Effectively, if the vertical component of an infinitesimal isometric deformation in Euclidean space
vanishes, then the deformation is trivial. (It comes from the fact that in
this case the rotation field associated to
the infinitesimal isometric deformation is constant, that is equivalent to the triviality of the deformation, see e.g.
\cite[Lemma 4, p. 256]{SPI5}). Proposition \ref{bordcon2} is now proved.

\subsubsection{Infinitesimal Pogorelov map}

The following construction is an adaptation of a map invented by Pogorelov 
\cite{Pogo}, which allows to transport deformation problems in 
a constant curvature space to deformation problems in a flat space, see for example
\cite{SchLab,rousset,Schconvex}.

We view hyperbolic space as a quadric in the Minkowski space of dimension $4$, that is
\begin{equation}
\mathbb{H}^3=\{ x\in\mathbb{R}^4 \vert x_1^2+x_2^2+x_3^2-x_4^2=-1, x_4>0\}.\nonumber
\end{equation}
We denote by $\varphi$  the projective map sending  $\mathbb{H}^3$  to 
the Klein projective  model. 
This map is known to be a homeomorphism between the hyperbolic
space and the unit open ball of the Euclidean space of dimension $3$, and it sends geodesics to straight
lines. Moreover, it sends convex sets to convex sets. In
particular, convex polyhedral surfaces are sent to convex polyhedral surfaces.

 Let $Z(x)$ be a vector of $T_x\mathbb{H}^3$. The \emph{radial component} of $Z(x)$ is the projection of $Z(x)$  on the radial direction, which is given by the derivative at $x$ of the geodesic $l_x$ in $\mathbb{H}^3$ between $x_c:=(0,0,0,1)$ (as a point of the Minkowski space) and $x$. The \emph{lateral component} of $Z(x)$ is the component orthogonal to the radial one. We denote by $\mu$  the length of the geodesic $l_x$, and by $Z_r$ and $Z_l$ the radial and lateral components of $Z$. The definitions are the same in Euclidean space, taking the origin instead of $x_c$ (the projective map sends $x_c$ to the origin).

The \emph{infinitesimal Pogorelov map} $\Phi$ is a map sending a vector field $Z$ of
hyperbolic space to a vector field $\Phi(Z)$ of Euclidean space, defined as follow:
the radial component of $\Phi(Z)(\varphi(x))$  has same direction and
same norm as $Z_r(x)$, and the lateral component of
$\Phi(Z)(\varphi(x))$ is $d_x\varphi (Z_l)$.

If we see a polyhedral surface $P$ in the Klein projective model, then the infinitesimal Pogorelov map is a map sending a vector field on $P$ to another vector field on $P$.

\begin{figure}[h] \begin{center}
\input{hyperboloide.pstex_t}
\end{center}
\caption{$\norm{Z_l}_{\mathbb{H}^3}=\cosh{\mu}\norm{\Phi(Z)_l}_{\mathbb{R}^3}$.\label{thales}}
\end{figure}
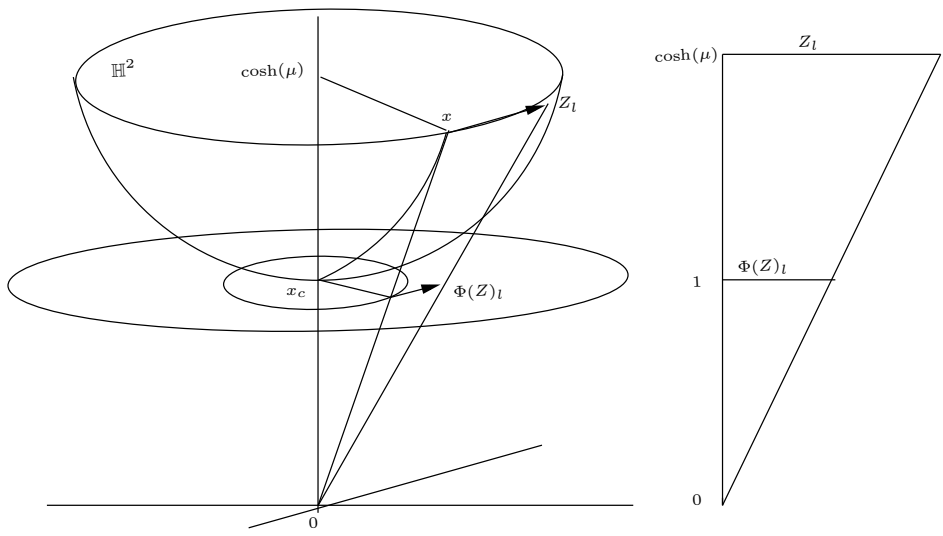

We have
\begin{eqnarray}\label{coef projection}
\ \norm{Z_r}_{\mathbb{H}^3}=\norm{\Phi(Z)_r}_{\mathbb{R}^3}; \
\norm{Z_l}_{\mathbb{H}^3}=\cosh{\mu}\norm{\Phi(Z)_l}_{\mathbb{R}^3}.\end{eqnarray}
The first one is the definition, the second one comes from a direct computation or an 
elementary property of the geometry of the plane (sometimes called the Thales Theorem, see Figure \ref{thales}).

We will sometimes make the confusion consisting to miss out the point at which we evaluate a vector field.

The infinitesimal Pogorelov map has the following remarkable property:
\begin{lemme}[Fundamental property of the infinitesimal Pogorelov map {\cite[1.9]{Schconvex}}]\label{pogorelov}
Let $V$ be 
a vector field on $\mathbb{H}^3$, then $V$ is a Killing field if and only if  $\Phi(V)$  is a Killing field of the Euclidean space.
\end{lemme}
As an infinitesimal isometric deformation of a polyhedral surface is the data of a Killing field on each triangle of a triangulation, this lemma says that
the image of an infinitesimal isometric deformation of a polyhedral surface $P$  by the infinitesimal Pogorelov map is an  infinitesimal isometric deformation of the image of $P$ by the projective map. And one is trivial when the other is.

\subsection{Fuchsian infinitesimal rigidity}

We want to show that convex  Fuchsian polyhedra are infinitesimally rigid among  Fuchsian polyhedra. As we
will consider the Fuchsian polyhedra up to isometries of hyperbolic space, we  consider that the
invariant totally geodesic surface in the definition  is always the same. We choose 
the horizontal plane in the Klein projective model, and we will denote it by
$P_{\mathbb{H}^2}$. By $p_{\mathbb{H}^2}$ we mean the orthogonal projection in hyperbolic space onto the
plane $P_{\mathbb{H}^2}$.

\paragraph{Global form of convex Fuchsian polyhedra.}

Let $(P,F)$ be a  convex Fuchsian polyhedron in the hyperbolic
space. 
\begin{lemme}
The polyhedral surface $P$ has a boundary at infinity which is the same as the one of $P_{\mathbb{H}^2}$.
\end{lemme}
\begin{proof}
It is a general property of discrete subgroups of isometries of hyperbolic space that, for a point $x\in\mathbb{H}^3$, the set of accumulation 
points of $(fx)_{f\in F}$ doesn't depend on the choice of the 
point $x$. As  $F$ acts on $P$ and on $P_{\mathbb{H}^2}$, we deduce from this that the boundary at infinity of $P$ is included 
in the boundary at infinity of $P_{\mathbb{H}^2}$.

It remains to check that the fact that $F$ acts cocompactly 
on  $P_{\mathbb{H}^2}$
implies that all the points of the boundary at infinity of $P_{\mathbb{H}^2}$ are reached by its action, and this is left to the reader.
\end{proof}

Then, by convexity of $P$, $P$ and $P_{\mathbb{H}^2}$ have no intersection point in hyperbolic space
(if we glue them along their common infinite boundary, they bound a convex body - recall that in the Klein projective model, $P_{\mathbb{H}^2}$ is the intersection of a Euclidean plane with the unit ball) and, again by  convexity of $P$, 
the orthogonal projection on $P_{\mathbb{H}^2}$ gives a homeomorphism between
$P_{\mathbb{H}^2}$ and $P$. Note that this implies that for a  convex Fuchsian polyhedron $(P,F)$, the genus of the surface 
$P/F$ is inevitably $>1$.

By cocompactness of the action of $F$ on $P_{\mathbb{H}^2}$, the distance between points of $P$ and points of  $P_{\mathbb{H}^2}$
(given by the orthogonal projection on $P_{\mathbb{H}^2}$)  is
bounded: $P$ is between two umbilic surfaces, realising
the extrema $d_{min}$ and $d_{max}$ of the distance between $P$ and
$P_{\mathbb{H}^2}$.

This leads to the fact that, if we see $P$ in the Klein projective model, it is a convex polyhedral
cap (with infinite number of vertices accumulating on
the  boundary), lying between two half-ellipsoids of radii
$(1,1,r_{min})$ and $(1,1,r_{max})$, with $0<r_{min}<r_{max}<1$ (a direct computation shows that in this model, umbilic surfaces 
are half-ellipsoids of radius $(1,1,r)$ where $r=tanh(d) <1$, with
$d$ the hyperbolic distance between the umbilic surface and $P_{\mathbb{H}^2}$).

\paragraph{Fuchsian polyhedral embeddings.}

We will need to introduce another way to describe Fuchsian polyhedra. 
In all the following, $S$ is a compact surface of genus $>1$. 
\begin{definition}
A \emph{polyhedral embedding} of $S$ in hyperbolic space is  a cellulation of $S$ together with a homeomorphism from 
$S$ to a polyhedral surface of the hyperbolic space $\mathbb{H}^3$, sending polygons of the cellulation
to geodesic polygons of  $\mathbb{H}^3$. 

A \emph{Fuchsian polyhedral  embedding}
of  $S$ in the hyperbolic space  $\mathbb{H}^3$ is a couple
 $(\phi,\rho)$, where
\begin{itemize}
\item $\phi$ is a polyhedral embedding
of the universal cover $\widetilde{S}$ of $S$ in $\mathbb{H}^3$,
\item  $\rho$ is a representation of the fundamental group $\Gamma$ of  $S$ in the group of orientation-preserving isometries of
$\mathbb{H}^ 3$,
\end{itemize}
such that $\phi$ is \emph{equivariant} under the action of $\Gamma:= \pi_1(S)$:
\begin{equation}
\ \forall \gamma \in \Gamma, \forall x \in \widetilde{S},
\: \phi(\gamma x)=\rho(\gamma)\phi(x),\nonumber
\end{equation}
and $\rho(\Gamma)$ leaves globally invariant a totally geodesic surface in $\mathbb{H}^3$, on which it acts cocompactly (without fixed points).

The \emph{number of vertices of the Fuchsian polyhedral  embedding} is the number of vertices of the cellulation 
of $S$.

The Fuchsian polyhedral embedding is \emph{convex} if its image is a convex polyhedral surface of the hyperbolic space.
\end{definition}

We consider the Fuchsian polyhedral  embeddings up to homeomorphisms and up to global isometries: let $(\phi_1,\rho_1)$ and $(\phi_2,\rho_2)$ be two  Fuchsian polyhedral embeddings of two surfaces $S_1$ and 
$S_2$. We say that  $(\phi_1,\rho_1)$ and $(\phi_2,\rho_2)$ are equivalent if there exists 
a homeomorphism $h$ between $S_1$ and $S_2$ and a hyperbolic isometry $I$
such that, for a lift $\widetilde{h}$ of $h$ to $\widetilde{S}_1$ we have
\begin{equation}
\phi_2\circ \widetilde{h} = I\circ \phi_1.\nonumber
\end{equation}
As two lifts of $h$ only differ by conjugation by elements of
$\Gamma:=\pi_1(S)$, using the equivariance property of the
embedding, it is easy to check that the definition of the
equivalence relation doesn't depend on the choice of the lift.

As we see the Fuchsian polyhedral  embeddings up to global isometries, we consider that the invariant surface is 
always $P_{\mathbb{H}^2}$.
\begin{definition}
The \emph{genus of a Fuchsian group $F$ of hyperbolic space}  is the genus 
of the quotient of the invariant totally geodesic surface (for the action of $F$) by the restriction of $F$ to it.

The \emph{genus of a Fuchsian polyhedron} is the genus of the Fuchsian group of hyperbolic space acting on it.

The 
\emph{number of vertices of a Fuchsian polyhedron} $(P,F)$ is the
number of vertices of $P$ in a fundamental domain for the action of 
$F$.
\end{definition}
As $S$ is a compact surface of genus $g>1$, it can be endowed with hyperbolic metrics, and each of them provides 
a cocompact representation of $\Gamma$ in the group of orientation-preserving  isometries of the hyperbolic plane. The images of such representations are usually called \emph{Fuchsian groups} (of $\mathbb{H}^2$), that explains the  terminology used. Moreover, there
is a bijection between the cocompact representations of $\Gamma$ 
in $\mbox{Isom}^+(\mathbb{H}^2)$ and the Fuchsian groups of $\mathbb{H}^2$ of genus $g$.

\begin{lemme}
There is a bijection between the cocompact representations of the fundamental group of $S$ in 
$\mbox{Isom}^+(\mathbb{H}^2)$ and 
the Fuchsian groups of $\mathbb{H}^3$ of genus $g$ (which leave invariant $P_{\mathbb{H}^2}$).
\end{lemme}
\begin{proof}
It suffices to prove that there is a bijection between the Fuchsian groups of $\mathbb{H}^2$ and the Fuchsian groups of $\mathbb{H}^3$ (with same genus).

The restriction of a Fuchsian group of $\mathbb{H}^3$ to $P_{\mathbb{H}^2}$ obviously
gives a Fuchsian group of $\mathbb{H}^2$.
  
  Reciprocally,  a Fuchsian group $F$  acting on $P_{\mathbb{H}^2}$ canonically gives a Fuchsian group 
  of the hyperbolic space: for a point $x$ in  the hyperbolic space,
an element $f\in F$ sends $p_{\mathbb{H}^2}(x)$ on a point $z\in P_{\mathbb{H}^2}$. The image $y$ of $x$ by the element of $\mbox{Isom}^+(\mathbb{H}^3)$ extending $f$ is the unique point of the 
hyperbolic space such that its projection on $P_{\mathbb{H}^2}$ is $z$ (and is in the same 
half-space delimited by $P_{\mathbb{H}^2}$ than $x$). And there is no other such subgroup of $\mbox{Isom}^+(\mathbb{H}^3)$ (because if there is, each elements of  both groups sends a geodesic segment orthogonal to 
$P_{\mathbb{H}^2}$ to the same geodesic segment, then they are equal).
\end{proof}

\begin{lemme}
There is a bijection between the convex Fuchsian polyhedra of genus $g$ with $n$ vertices and the 
convex Fuchsian polyhedral embeddings with $n$ vertices of $S$.
\end{lemme}
\begin{proof}

Obviously,  the image of $S$ by a convex Fuchsian polyhedral  embedding is a convex Fuchsian polyhedron.
Reciprocally, the canonical embedding in  $\mathbb{H}^3$ 
of a convex Fuchsian polyhedron $P$ invariant under the action 
of a group  $F$  
gives a convex Fuchsian polyhedral  embedding of the surface $P/F$ in $\mathbb{H}^3$. 
We have seen that this surface is homeomorphic to $P_{\mathbb{H}^2}/F$, which is homeomorphic to $S$.
\end{proof}

\paragraph{Fuchsian deformations.}

Let $(S,\phi,\rho)$ a   convex  polyhedral Fuchsian embedding.

Let $(\phi_t)_t$ be  a path of  convex polyhedral embeddings   of
$\widetilde{S}$ in  $\mathbb{H}^3$, such that:
\begin{itemize}
\item[-] $\phi_0=\phi$,
\item[-] the induced metric is preserved at the first order at $t=0$,
\item[-] there are representations $\rho_t$ of $\Gamma=\pi_1 (S)$ in $\mbox{Isom}^+(\mathbb{H}^3)$
\end{itemize}
 such that
\begin{equation}\phi_t(\gamma x)=\rho_t(\gamma)\phi_t(x)\nonumber\end{equation}
and each $\rho_t(\Gamma)$ leaves globally invariant  a totally geodesic surface, on which it acts cocompactly without fixed points
(up to global isometries, we consider that the surface is always  $P_{\mathbb{H}^2}$). 

We denote by \begin{equation}Z(\phi(x)):= \frac{d}{dt}\phi_t(x)_{\vert t=0}\in T_{\phi(x)}\mathbb{H}^3\nonumber\end{equation} and
\begin{equation}\dot{\rho}(\gamma)(\phi(x))=\frac{d}{dt}\rho_t(\gamma)(\phi(x))_{\vert t=0} \in T_{\rho(\gamma)\phi(x)}\mathbb{H}^3.\nonumber\end{equation}
The vector field $Z$ has a property of equivariance under
$\rho(\Gamma)$:
\begin{equation}\label{equ definition def}
Z(\rho(\gamma) \phi(x))=\dot{\rho}(\gamma)(\phi(x))+d\rho(\gamma).Z(\phi(x)).
\end{equation}
This can be written
\begin{equation}\label{killequiv}
Z(\rho(\gamma)\phi( x))=d\rho(\gamma).(d\rho(\gamma)^{-1}\dot{\rho}(\gamma)(\phi(x))+Z(\phi(x)))
\end{equation}
and $d\rho(\gamma)^{-1}\dot{\rho}(\gamma)$ is a Killing field of  $\mathbb{H}^{3}$, because it is the derivative of a path in the group of isometries of  $\mathbb{H}^{3}$ (we must multiply by  $d\rho(\gamma)^{-1}$, because
$\dot{\rho}(\gamma)$ is not a vector field). We denote this Killing field by $\vec{\rho}(\gamma)$. Equation (\ref{killequiv}) can be written, if $y=\phi(x)$,
\begin{equation}\label{equation deformation fuchsienne}
Z(\rho(\gamma) y)=d\rho(\gamma).(\vec{\rho}(\gamma) + Z)(y).
\end{equation}
A \emph{Fuchsian deformation} is an infinitesimal isometric deformation  $Z$ on a Fuchsian polyhedron which satisfies Equation (\ref{equation deformation fuchsienne}), where $\vec{\rho}(\gamma)$ is a \emph{Fuchsian Killing field}, that is a Killing field of the hyperbolic plane extended to the hyperbolic space along geodesics orthogonal to the plane. More precisely, for a point  $x\in\mathbb{H}^3$, let  $d$ be the
distance between $x$ and $p_{\mathbb{H}^2}(x)$. We denote by $p_d$
the orthogonal projection onto $P_{\mathbb{H}^2}$ of the umbilic
surface which is at constant distance $d$ from $P_{\mathbb{H}^2}$ (passing
through $x$). Then the Killing field $K$ at $p_{\mathbb{H}^2}(x)$ 
is extended as $dp_d^{-1}(K)$ at the point $x$.

A Fuchsian polyhedron  is \emph{Fuchsian infinitesimally  rigid}  if all its Fuchsian deformations are trivial (i.e. are restriction to the Fuchsian polyhedron of Killing fields of hyperbolic space).

We want to prove
\begin{theoreme}\label{theoreme de rigidite}
Convex Fuchsian polyhedra are Fuchsian infinitesimally  rigid.
\end{theoreme}

By the fundamental Property of the infinitesimal Pogorelov map, 
to prove Theorem \ref{theoreme de rigidite} it suffices to prove
that, for each Fuchsian deformation $Z$, the vertical component of
the image of $Z$ by the infinitesimal Pogorelov map vanishes at the
boundary, because in this case Proposition \ref{bordcon2} provides the result.

\paragraph{Remarks about the method employed.}

There is a non-direct  way to prove the statement of Theorem \ref{theoreme de rigidite}, 
using another infinitesimal Pogorelov map. This way is pointed out in \cite{artrealisationlorentz}.

Moreover, it may be possible that the method employed to prove  Proposition \ref{bordcon2} 
leads to a direct proof of Theorem \ref{theoreme de rigidite} in hyperbolic space, without using the infinitesimal Pogorelov map.

But we think that the method used here can be extended to prove infinitesimal rigidity results for ideal and hyperideal Fuchsian polyhedra (that means that some vertices could be on the sphere or out of the ball 
for the Klein projective model).

Note that the following proof is also true without any change for strictly convex smooth Fuchsian surfaces (using the smooth analog of Proposition \ref{bordcon2}).

\paragraph{Proof of Theorem \ref{theoreme de rigidite}.}
Recall that $p_{\mathbb{H}^2}$ is the orthogonal projection  on the totally geodesic surface $P_{\mathbb{H}^2}$.  At a point $x$ 
of a convex Fuchsian polyhedron $P$, the derivative of the geodesic which realises this projection is called the \emph{vertical direction} at $x$, and the directions orthogonal to this one are \emph{horizontal directions}. So a Fuchsian deformation
$Z$ can be decomposed in a vertical component $Z_v$ and in a horizontal component $Z_h$.
We denote by $(Z_r)_h$ the horizontal component of the radial component of $Z$, etc.
We have
\begin{equation}
Z_r=(Z_r)_h+(Z_r)_v=(Z_h)_r+(Z_v)_r.\nonumber
\end{equation}
The first one is obvious and the second one comes from the linearity of the orthogonal projection.

\begin{proposition}
The vector field  $Z_v$ is \emph{invariant} under the action of $\rho(\Gamma)$, i.e. $\forall x\in P$:
\begin{equation}Z_v(\rho(\gamma) x)=d\rho (\gamma). Z_v(x).\nonumber\end{equation}
\end{proposition}
\begin{proof}
From Equation (\ref{equation deformation fuchsienne}), it suffices to check that  $\vec{\rho}(\gamma)$ has no vertical component, that is true by definition of a Fuchsian Killing field.
\end{proof}
And as  $\rho(\Gamma)$ acts cocompactly on $P$
\begin{corollaire}\label{la composante verticale est bornee}
There exists a constant  $c_v$ such that, for all $x\in P$,
\begin{equation}\norm{Z_v(x)}\leqslant c_v.\nonumber
\end{equation}
\end{corollaire}
Moreover
\begin{corollaire}
The vector field $Z_h$ is equivariant under the action of $\rho(\Gamma)$.
\end{corollaire}
\begin{proof}
 Equation (\ref{equ definition def}) says that
\begin{equation}Z_v(\phi(\gamma x))+Z_h(\phi(\gamma x))=\dot{\rho}(\phi (x)) +
 d\rho (\gamma).Z_v(\phi(x))+d\rho(\gamma).Z_h(\phi(x)), \nonumber\end{equation}and the preceding proposition
gives the result.
\end{proof}

Recall that $p_d$ is
the orthogonal projection on $P_{\mathbb{H}^2}$ of the umbilic
surface which is at constant distance $d$ from $P_{\mathbb{H}^2}$ (passing
through $x$). We call \emph{radial-horizontal} the component of
$Z$ (at $x$) in the direction $dp_d^{-1}(r(p_{\mathbb{H}^2}(x)))$, where $r(p_{\mathbb{H}^2}(x))$ is
the radial direction of $P_{\mathbb{H}^2}$ at the point
$p_{\mathbb{H}^2}(x)$. This component is noted $Z_{rh}$, and it's a horizontal vector.

We denote by $W$ the projection on $P_{\mathbb{H}^2}$ of the horizontal component of $Z$ (it is equivariant under the action of $\rho(\Gamma)$).
We denote by $W_{r}$ its radial component.
Then $dp_d^{-1}(W_r)$ is the radial-horizontal component of the horizontal component of  $Z$.

The determining fact is:
\begin{proposition}\label{majproj}
Let $H$ be a vector field of  $\mathbb{H}^{2}$ equivariant under the action of  $\rho(\Gamma)$. Then there exists
a constant  $c_{\dot{h}}$ such that
\begin{equation}\norm{H_{r}(x)}_{\mathbb{H}^{2}}\leq c_{\dot{h}} d_{\mathbb{H}^{2}}(x_{c},x).\nonumber\end{equation}
\end{proposition}
The point $x_c$ is always the origin in the Klein projective model.
\begin{proof}
We deform the hyperbolic metric $h$ along  $H$, that is  $\dot{h}=L_{H}h$, where $L$ is the Lie derivative:
\begin{eqnarray}
 \ \dot{h}(X,Y)=(L_{H}h)(X,Y)&=&h(\nabla_{X}H,Y)+h(\nabla_{Y}H,X)\nonumber \\
  \ &=&H.h(X,Y)+h([X,H],Y)+h(X,[Y,H]).\nonumber
\end{eqnarray}
Let  $c\: : \: [0,\eta]\rightarrow \mathbb{H}^{2},\: c(0)=x_{c},\:, c(\eta)=x,\:\norm{c'}=1$ a geodesic (then $\eta=d_{\mathbb{H}^2}(x_c,x)$ and $c'(\eta)$ is the radial direction 
at the point $x$).
Up to adding a Killing field,
we can suppose that  $H(x_c)=0$. Then

\begin{eqnarray}
\ \int_{0}^{\eta}\dot{h}(c'(t),c'(t))dt&=&2\int_{0}^{\eta}h(\nabla_{c'(t)}H(c(t)),c'(t))dt \nonumber \\
\ &=&2\int_{0}^{\eta}c'(t).h(H,c'(t))dt \nonumber \\
\ &=&  2h(H,c'(\eta))\nonumber \\
\ &=& 2h(H_r+H_l,c'(\eta)) \nonumber \\
\ &=& 2 h(H_r,c'(\eta))\nonumber \\
\ &=& 2 h(H_r,\pm\frac{H_r}{\norm{H_r}_{\mathbb{H}^2}}),\nonumber 
\end{eqnarray}
that leads to
\begin{equation}\label{majov}
\int_{0}^{\eta}\dot{h}(c'(t),c'(t))dt=\pm 2 \norm{H_{r}}_{\mathbb{H}^{2}}.
\end{equation}
As  $H$ is equivariant under the action of  $\rho(\Gamma)$,  the elements of  $\rho(\Gamma)$ preserve the bilinear form $\dot{h}$, because,
up to an isometry, $H(\rho(\gamma)x)$ is written as  $H(x)$ plus a Killing field.

Formally, using Equation (\ref{equation deformation fuchsienne})
\begin{equation}
H(\rho(\gamma)x)=d\rho(\gamma).(\vec{\rho}(\gamma)+H)(x)\nonumber,
\end{equation}
and we develop
\begin{eqnarray}
 \ \dot{h}(d\rho(\gamma)X(x),d\rho(\gamma)Y(x)) &:=&h(\nabla_{d\rho(\gamma)X(x)}H(\rho(\gamma)x),d\rho(\gamma)Y(x)) \nonumber\\
\ &&+  h(\nabla_{d\rho(\gamma)Y(x)}H(\rho(\gamma)x),d\rho(\gamma)X(x))  \nonumber\\
\ &=& h(\nabla_{d\rho(\gamma)X(x)}  d\rho(\gamma) ( \vec{\rho}(\gamma)(x)+H(x)),d\rho(\gamma)Y(x))  \nonumber\\
\ &&+ h(\nabla_{d\rho(\gamma)Y(x)} d\rho (\gamma) (\vec{\rho}(\gamma)(x)+H(x)),d\rho(\gamma)X(x))  \nonumber\\
\ &=&h(\nabla_{X(x)}  ( \vec{\rho}(\gamma)(x)+H(x)),Y(x))  \nonumber\\
\ &&+ h(\nabla_{Y(x)}  (\vec{\rho}(\gamma)(x)+H(x)),X(x))  \nonumber\\
\ &=& h(\nabla_{X(x)}H(x),Y(x))+ h(\nabla_{Y(x)}H(x),X(x)) \nonumber\\
\ &&+h(\nabla_{X(x)}\vec{\rho}(\gamma)(x),Y(x))+ h(\nabla_{Y(x)}\vec{\rho}(\gamma)(x),X(x))  \nonumber\\
\ &=& h(\nabla_{X(x)}H(x),Y(x))+ h(\nabla_{Y(x)}H(x),X(x)) \nonumber\\
\ &&+L_{\vec{\rho}(\gamma)}h(X(x),Y(x))  \nonumber\\
\ &=& \dot{h}(X(x),Y(x)) \nonumber
\end{eqnarray}
($L_{\vec{\rho}(\gamma)}h(X,Y)=0$ because $\vec{\rho}(\gamma)$ is a Killing field).

Then $\dot{h}$ is a bilinear form on  $P_{\mathbb{H}^2}/\rho(\Gamma)$, which is compact, then $\vert \dot{h}\vert$ is bounded by a constant
$c_{\dot{h}}$, and by Equation (\ref{majov}):
\begin{equation}
\norm{H_{r}(x)}_{\mathbb{H}^{2}}\leq c_{\dot{h}} \eta =c_{\dot{h}} d_{\mathbb{H}^{2}}(x_c,x).\nonumber
\end{equation}
\end{proof}

\begin{corollaire}\label{trucH2}
There exists a constant $c_{rh}$ such that:
\begin{equation}\norm{(Z_{h})_{rh}(x)}_{\mathbb{H}^3}\leq c_{rh} d_{\mathbb{H}^3}(x_c,x).\nonumber
\end{equation}
 \end{corollaire}
\begin{proof}
A simple computation in Minkowski space shows that the induced metric on the umbilic surface at constant distance $d$ of $P_{\mathbb{H}^2}$ (passing through $x$) is
$\cosh (d) h$ where  $h$ is the hyperbolic metric, then
\begin{equation}
\norm{(Z_h)_{rh}(x)}_{\mathbb{H}^3}= \cosh (d) \norm{W_{r}(p_{\mathbb{H}^2}(x))}_{\mathbb{H}^2},\nonumber
\end{equation}
and as $d$ is bounded for all $x\in P$, and with the preceding proposition, there exists a constant  $c_{rh}$ such that:
\begin{equation}\label{equation majoration intermediaire}
\norm{(Z_h)_{rh}(x)}_{\mathbb{H}^3}\leq c_{rh} d_{\mathbb{H}^2}(x_{c},p_{\mathbb{H}^2}(x)).
\end{equation}

Classical hyperbolic trigonometry %(\cite[2.6]{Thurcour1}), 
applied to the rectangular triangle with edge length  $d_{\mathbb{H}^3}(x_{c},x)$ (the long edge),    $d_{\mathbb{H}^3}(x,p_{\mathbb{H}^2}(x))$ and  $d_{\mathbb{H}^2}(x_c,p_{\mathbb{H}^2}(x))$ gives:
\begin{equation}\cosh (d_{\mathbb{H}^3}(x_c,x))=\cosh( d_{\mathbb{H}^3}(x,p_{\mathbb{H}^2}(x))) \cosh (d_{\mathbb{H}^2}(x_c,p_{\mathbb{H}^2}(x))),\nonumber\end{equation}
and as the values of the hyperbolic cosine are greater than $1$:
\begin{equation}\cosh (d_{\mathbb{H}^3}(x_{c},x)) \geq \cosh (d_{\mathbb{H}^2}(x_{c},p_{\mathbb{H}^2}(x))),\nonumber\end{equation}
and the hyperbolic  cosine is an increasing function for positive values, then:
\begin{equation}d_{\mathbb{H}^3}(x_{c},x) \geq d_{\mathbb{H}^2}(x_{c},p_{\mathbb{H}^2}(x)).\nonumber\end{equation}
This together with Equation (\ref{equation majoration intermediaire}) gives the result.
\end{proof}

The  \emph{vertical plane} is the vector space of  $T_{x}\mathbb{H}^{3}$ spanned by the orthogonal vectors $Z_{v}$
and $Z_{rh}$ ($Z_{rh}$ is in the horizontal plane by definition, then it is orthogonal to the vertical direction).

We can see the vertical plane as the tangent plane  (at $x$) to the totally geodesic surface passing through  $x_c$ and $x$ and orthogonal
 to $P_{\mathbb{H}^2}$.

We denote by $Z_{lv}$ the image of the projection on the vertical plane of the
lateral component of  $Z$.

The vector $Z_r$ belongs to the vertical plane, because it can be decomposed in a horizontal component, which is in the radial-horizontal direction, and a vertical component.

The lateral component is orthogonal to the radial component, then
the vector $Z_{lv}$ is orthogonal to  $Z_{r}$ in the vertical
plane.

In all the following, we denote by $\mu:=d_{\mathbb{H}^{3}}(x_c,x)$,
and vector fields are evaluated at the point $x\in P$.

\begin{lemme}\label{decomposition vecteur vertical}
Let $V$ be the projection of a component of $Z$ on the vertical plane. Then there exists a constant  $c$ such that
\begin{equation}
\norm{V}_{\mathbb{H}^3} \leq c (1+\mu).\nonumber
\end{equation}
\end{lemme}
\begin{proof}
We denote by $\Pi_{V}$ the projection on the vertical plane, considered as spanned by the orthogonal vectors  $Z_{rh}$ and $Z_{v}$.
 We can write $\Pi_{V}(Z)=Z_{rh}+Z_{v}$. As $V$ is already in the vertical plane,
and as we project a component of  $Z$, we can write:
 \begin{eqnarray*}
 \ \norm{V}_{\mathbb{H}^{3}}=\norm{\Pi_{V}(V)}_{\mathbb{H}^{3}}
 &\leqslant& \norm{\Pi_{V}(Z)}_{\mathbb{H}^{3}}  \\
  \ &\leqslant&  \norm{Z_{rh}}_{\mathbb{H}^{3}}+ \norm{Z_{v}}_{\mathbb{H}^{3}} \\
 \ &\leqslant&  \norm{(Z_{h})_{rh}}_{\mathbb{H}^{3}}+
 \norm{(Z_{v})_{rh}}_{\mathbb{H}^{3}}+ \norm{Z_{v}}_{\mathbb{H}^{3}} \\
  \ &\leqslant&  \norm{(Z_{h})_{rh}}_{\mathbb{H}^{3}}+
  2\norm{Z_{v}}_{\mathbb{H}^{3}},
   \end{eqnarray*}
 and as the overestimation of these two last norms are known (by Corollaries \ref{la composante verticale est bornee} and
 \ref{trucH2}) we get
\begin{equation}
\norm{V}_{\mathbb{H}^3} \leqslant c_{rh} \mu +2c_v,\nonumber
\end{equation}
that is, if  $c$ is greater than  $c_{rh}$ and $2c_{v}$,
\begin{equation}
\norm{V}_{\mathbb{H}^3} \leqslant c(1+ \mu).\nonumber
\end{equation}
\end{proof}

For convenience, we denote by  $u$ be the image of $Z$ by the infinitesimal Pogorelov
map.

The decompositions of tangent vectors defined above (vertical, horizontal,\ldots) occur in Euclidean space. We want to show that the vertical component of $u$ goes to $0$ at the boundary.
 As $u_{v}$ is in the vertical plane, if  $\alpha$ is the angle between $u_{v}$ and $u_{lv}$, we get
\begin{equation}\label{ecrituresin}
u_{v}=\cos(\alpha) u_{lv}+ \sin({\alpha}) u_{r},
\end{equation}
because, for the same reason as  for $Z$,  $u_{r}$ and $u_{lv}$ give an
orthogonal basis of the vertical plane.

By  Lemma \ref{decomposition vecteur vertical} we get, as $Z_r$ and $Z_{lv}$ are in the vertical plane,
\begin{eqnarray*}
\ &&\norm{Z_r}_{\mathbb{H}^3} \leqslant c(1+\mu), \\
\  &&\norm{Z_{lv}}_{\mathbb{H}^3} \leqslant c(1+\mu),
\end{eqnarray*}
and as the infinitesimal Pogorelov map preserves the norm of the radial
component and crushes by
a coefficient
 $\cosh(\mu)$  the norm of the lateral direction (Equation (\ref{coef projection})), these two inequalities rise to
\begin{eqnarray*}
\ &&\norm{u_r}_{\mathbb{R}^3} \leqslant c(1+ \mu), \\
\ &&\norm{{u}_{lv}}_{\mathbb{R}^3} \leqslant \cosh(\mu)^{-1} c(1+\mu).
\end{eqnarray*}

Starting from (\ref{ecrituresin})
we get
\begin{eqnarray*}
\ \norm{u_{v}}_{\mathbb{R}^{3}}&\leqslant& \norm{u_{lv}}_{\mathbb{R}^{3}}
+ \sin({\alpha}) \norm{u_{r}}_{\mathbb{R}^{3}} \\
\ &\leqslant& c\cosh(\mu)^{-1}(1+ \mu) + c\sin({\alpha}) (1+ \mu)
 \\
 \ &\leqslant& c(1+ \mu) (\cosh(\mu)^{-1}+\sin({\alpha})).
\end{eqnarray*}

We call $\delta$ the Euclidean distance between a point $x$ on $P$
and $P_{\mathbb{H}^2}$ (in the Klein projective model), and we see easily that we have
  $\sin(\alpha)\underset{0}\thickapprox\delta$ when we go near the boundary of the ball:
\begin{equation}
c(1+ \mu) (\cosh(\mu)^{-1}+\sin({\alpha}))\underset{0}\thickapprox c(1+ \mu) (\cosh(\mu)^{-1}+\delta).\nonumber
\end{equation}

We denote by $S_{max}$  the umbilic surface which realises the
maximum of the distance between $P$ and $P_{\mathbb{H}^2}$,
$x_{max}$ the intersection of $S_{max}$ with the geodesic joining
$x$ and $x_c$, $\mu_{max}$ for
$d_{\mathbb{H}^3}(x_c,x_{max})$, $\delta_{max}$ the distance in
$\mathbb{R}^{3}$ between  $x_{max}$ and $P_{\mathbb{H}^2}$.
Guess analog definitions of $S_{min}, x_{min},\mu_{min}, \delta_{min}$
for the surface realising the minimum of the distance between
$P$ and $P_{\mathbb{H}^2}$.

\begin{lemme}\label{lemmeapprox}
Near the boundary of the ball, we have the approximations
\begin{eqnarray*}
\  \mu_{max} \underset{\delta_{max}
 \rightarrow 0}{\thickapprox}- c_{max}\ln(\delta_{max}), \\
\ \mu_{min} \underset{\delta_{min}
 \rightarrow 0}{\thickapprox}- c_{min}\ln(\delta_{min}),
\end{eqnarray*}
 where $c_{max}$ and $c_{min}$ are positive constants.
\end{lemme}
\begin{proof}
We prove the lemma in the case $max$, the proof for  $min$ is the same.
It is easy to check that  $\mu_{max}=d_{\mathbb{H}^3}(x_c,x_{max})=\tanh^{-1}(\norm{x_{max}}_{\mathbb{R}^{3}})$, that is
\begin{equation}\mu_{max}= \ln \left( \frac{1+\norm{x_{max}}_{\mathbb{R}^3}}{1-\norm{x_{max}}_{\mathbb{R}^3}}\right)
\underset{\norm{x_{max}}_{\mathbb{R}^3}
 \rightarrow 1}{\thickapprox} - \ln(1-\norm{x_{max}}_{\mathbb{R}^{3}}).\nonumber
 \end{equation}
As the image of  $S_{max}$ in the projective Klein model is an ellipsoid,
$\delta_{max}$ satisfies the equation
\begin{equation}
(x_{max})_{1}^{2}+(x_{max})_{2}^{2}+\frac{\delta_{max}^{2}}{r^{2}}=1,\nonumber\end{equation}
where $r$ is a positive constant strictly less than  $1$.
 Adding and removing a  $\delta_{max}^{2}$ and reordering we get:
  \begin{equation}
  1-\norm{x_{max}}_{\mathbb{R}^{3}}=\delta_{max}^{2}\frac{1-r^{2}}{r^{2}},\nonumber
  \end{equation}
as $\delta_{max}^2=(x_{max})_3^2$,
and this gives the result.
 \end{proof}

As obviously
\begin{equation}
\mu_{min} \leq \mu\leq \mu_{max} \mbox{ and } \delta_{min} \leq \delta \leq \delta_{max},\nonumber
\end{equation}
we get
\begin{eqnarray*}
 \ \norm{u_v}_{\mathbb{R}^3} &\leq & c(1+ \mu) (\cosh(\mu)^{-1}+\delta) \\
\ &\leq &  c(1+ \mu_{max}) (\cosh(\mu_{min})^{-1}+\delta),
\end{eqnarray*}
and when $x$  goes near the boundary,   $\cosh(f(x))\underset{\infty}{\thickapprox} \exp(f(x))$, where $f$ is a function going to
$\infty$ when   $x$ goes near the boundary. Then:
\begin{eqnarray*}
\   && c(1+ \mu_{max}) (\cosh(\mu_{min})^{-1}+\delta) \\
\  \underset{0}\thickapprox &&  c(1- c_{max}\ln(\delta_{max})) (\cosh(- c_{min}\ln(\delta_{min}))^{-1}+\delta) \\
\ \underset{0}\thickapprox &&c'(1- \ln(\delta_{max})) (\delta_{min}+\delta).
\end{eqnarray*}
At the end, as
\begin{equation}
c'(1- \ln(\delta_{max})) (\delta_{min}+\delta) \leq c'(1- \ln(\delta_{max})) (\delta_{max}+\delta_{max})\nonumber
\end{equation}
and as when $\delta$ goes to $0$, $\delta_{max}$ goes to $0$, then
    $\norm{u_v}_{\mathbb{R}^3}$ goes to $0$.

    Theorem \ref{theoreme de rigidite} is now proved.

%%%%%%%%%%%%%%%%%%%%%%%%%%%%%%
%%%%%%%%%%%%%%%%%%%%%%%%%%%%%%
%%%%%%%%%%%%%%%%%%%%%%%%%%%%%%%%
\section{Realisation of metrics}

\subsection{Set of Fuchsian polyhedra}

We denote by $\mathcal{P}(n)$ the set of  convex Fuchsian
polyhedral embeddings with $n$ vertices of a compact surface $S$  in $\mathbb{H}^3$,
modulo isotopies of  $S$ fixing the vertices of the
cellulation and modulo the isometries of $\mathbb{H}^3$.

More precisely, the equivalence relation is the following: let $(\phi_1,\rho_1)$ and $(\phi_2,\rho_2)$ be two convex Fuchsian polyhedral embeddings of $S$. We say that  $(\phi_1,\rho_1)$ and $(\phi_2,\rho_2)$ are equivalent if there exists 
\begin{itemize}
\item a homeomorphism $h$ of $S$ isotopic to the identity, such that if $h_t$ is the isotopy (i.e. $t\in[0,1]$, 
$h_0=h$ and  $h_1=id$), then $h_t$ fixes the vertices of the cellulation for all $t$,
\item  a hyperbolic isometry $I$,
 \end{itemize}
 such that, for
a lift $\widetilde{h}$ of $h$ to $\widetilde{S}$ we have
\begin{equation}
\phi_2\circ \widetilde{h} = I\circ \phi_1.\nonumber
\end{equation}
Here again, the definition of the
equivalence relation doesn't depend on the choice of the lift.

\paragraph{Z-V-C  coordinates for Teichm\"uller space.}

For more details about  Z-V-C coordinates (Z-V-C  stands for Zieschang--Vogt--Coldewey, \cite{ZVC}) we refer to \cite[6.7]{Buser}.
\begin{definition}
Let $g\geq 2$. A (geodesically convex) polygon of the hyperbolic plane with edges  (in the direct order)  $b_1,b_2,\overline{b}_1,\overline{b}_2,b_3,b_4,\ldots,\overline{b}_{2g}$ and with interior angles $\theta_1,\overline{\theta}_1,\ldots,\theta_{2g},\overline{\theta}_{2g}$ is called   \emph{(normal) canonical} if, with $l(c)$ the length of the geodesic $c$,
\begin{itemize}
\item[i)] $l(b_k)=l(\overline{b}_k)$, $\forall k$;
\item[ii)] $\theta_1+\ldots+\overline{\theta}_{2g}=2\pi$;
\item[iii)] $\theta_1+\theta_2=\overline{\theta}_1+\theta_2=\pi$.
\end{itemize}
Two canonical polygons $P$ and $P'$ with edges $b_1,\ldots,\overline{b}_{2g}$ and  $b'_1,\ldots,\overline{b}'_{2g}$ are said
\emph{equivalent} if there exists an isometry from $P$ to $P'$ such that the edge  $b_1$ is sent to the edge $b'_1$ and $b_2$ is sent to $b'_2$.
\end{definition}

If we identify the edges $b_i$ with the edges $\overline{b}_i$, we
get a compact  hyperbolic surface of genus $g$. This surface could
also be written $\mathbb{H}^2 / F$, where $F$ is the sub-group of
$\mbox{PSL}(2,\mathbb{R})=\mbox{Isom}^+(\mathbb{H}^2)$ generated by the translations along the edges
$b_i$ (the translation length is the length of $b_i$).
The interior of the polygon is a fundamental domain for the action
of $F$. This leads to a description of the Teichm\"uller space
$T_g$:
\begin{proposition}[{\cite[6.7.7]{Buser}}]
Let $\mathcal{P}_g$  be the set of equivalence classes of canonical polygons. An element of  $\mathcal{P}_g$
is described by the $(6g-6)$ real numbers (the Z-V-C coordinates):
\begin{equation}
(b_3,\ldots,b_{2g},\theta_3,\overline{\theta}_3,\ldots,\theta_{2g}, \overline{\theta}_{2g}).\nonumber
\end{equation}
Endowed with this topology, $\mathcal{P}_g$ is in analytic bijection with $T_g$.
\end{proposition}

\paragraph{Surjection on the Teichm\"uller space with marked points.}

The Teichm\"uller space of $F_{g,n}$ (i.e. of a compact surface of genus $g>1$ with $n>0$ marked points), denoted by $T_g(n)$, can be defined as the set of hyperbolic metrics on $F_{g,n}$  modulo isotopies of $F_{g,n}$, such that each isotopy $h_t$ fixes the marked points for all $t$.

Let $(\phi_1,\rho_1)$ and $(\phi_2,\rho_2)$ be two equivalent convex Fuchsian polyhedral embeddings of $S$ with $n$ vertices. Recall that
$h$ is an isotopy of $S$ and $\widetilde{h}$ its lift to $\widetilde{S}$. As  $h$ is homotopic to the identity,  $\forall x\in\widetilde{S}, \forall \gamma\in\Gamma$, we get:
\begin{eqnarray*}
\ &&\phi_2(\widetilde{h}(\gamma x))=I(\phi_1 \gamma x) \\
\ \Leftrightarrow && \phi_2((h)_*(\gamma)x)=I(\rho_1(\gamma)(\phi_1 (x)) \\
\ \Leftrightarrow && \phi_2(\gamma x)=I(\rho_1(\gamma)(\phi_1 (x)) \\
\ \Leftrightarrow && \rho_2(\gamma)(\phi_2(x))=I(\rho_1(\gamma)(\phi_1 (x)) \\
\ \Leftrightarrow && \rho_2(\gamma)(I(\phi_1(x)))=I(\rho_1(\gamma)(\phi_1 (x)).
\end{eqnarray*}

But if two orientation-preserving isometries of the hyperbolic
space are equal on an open set of a totally geodesic surface
(a face of the Fuchsian polyhedron), they are equal, then for
all $\gamma\in\Gamma$, $\rho_2(\gamma)=I\circ \rho_1
(\gamma)\circ I^{-1}$. As $\rho_1$ and $\rho_2$ are also
representations of  $\Gamma$ in $\mbox{PSL}(2,\mathbb{R})$ (modulo
conjugation by an element of $\mbox{PSL}(2,\mathbb{C})=\mbox{Isom}^+(\mathbb{H}^3)$), we deduce that
$\rho_1$ and $\rho_2$ are the same element of $\mbox{Hom}(\Gamma,
\mbox{PSL}(2,\mathbb{R}))/\mbox{PSL}(2,\mathbb{R})$.

Then in the hyperbolic plane $P_{\mathbb{H}^2}$, the canonical polygons associated to these two representations 
are equal, and, up to an isometry, the projection of the vertices of the Fuchsian polyhedra $\phi_1(\tilde{S})$ and $\phi_2(\tilde{S})$ on $P_{\mathbb{H}^2}$ gives the same $n$ marked points in this canonical polygon: we have described a map  $\mathcal{S}$  which to each element of  $\mathcal{P}(n)$ associates an element of $T_g(n)$.

And as we have seen that from any Fuchsian group and any
$n$ points on the plane, we can build a convex Fuchsian
polyhedron with $n$ vertices (it is enough  to take  $n$ points at
same distance from the plane),  $\mathcal{S}$  is surjective.

\paragraph{Manifold structure on  $\mathcal{P}(n)$.}

\begin{lemme}\label{ouvert ensemble sommet}
Let $[h]\in T_g(n)$. Then $\mathcal{S}^{-1}([h])$ is diffeomorphic to  the open unit ball of $\mathbb{R}^{n}$.
\end{lemme}
\begin{proof}

We will show that $\mathcal{S}^{-1}([h])$ is a contractible open subset of  $(\mathbb{R}_+)^{n}$.

We fix an element of $T_g(n)$, that is the action of a cocompact Fuchsian group $F$ on 
$P_{\mathbb{H}^2}$ and $n$ points $(y_1,\dots,y_n)$.
We denote by $(x_1,\ldots,x_n)$ $n$ points of the hyperbolic space such that their projection
on the plane $P_{\mathbb{H}^2}$ is exactly $(y_1,\dots,y_n)$. We look at the Klein projective model of the hyperbolic space, and the coordinates
of the points are those of $\mathbb{R}^3$. The orthogonal  projection, Euclidean or hyperbolic, of
a point on the horizontal plane gives the same point.

We want the vertices of the polyhedron obtained
as the boundary of the closure of the convex hull of the points $f x_i$, for any $f\in F$ and $i=1\ldots n$,
to be exactly the set of points of the form $f x_i$, for any $f\in F$ and $i=1\ldots n$. This is the same as saying that no point of the form $f x_i$ is in the interior of the convex hull of the others points of this form; that means, for all  $f x$ ($x\in\{x_1,\ldots,x_n\}$), and for all $f_ix_i,
f_jx_j,f_kx_k $  such that $p_{\mathbb{H}^2}(f x)$
is contained in the triangle formed by the points
$p_{\mathbb{H}^2}(f_i x_i)$, $p_{\mathbb{H}^2}(f_j x_j)$
and $p_{\mathbb{H}^2}(f_k x_k)$, then the plane generated by   $f_ix_i, f_jx_j,f_kx_k $ (strictly)
separates   $f x$ from $P_{\mathbb{H}^2}$.

Actually, it suffices to verify this condition only for the points  $(x_1,\ldots,x_n)$, because if the  plane generated by  $$f_ix_i, f_jx_j,f_kx_k $$ separates $f x$ from $P_{\mathbb{H}^2}$, then the  plane generated by    $$f^{-1}f_ix_i,f^{-1}
f_jx_j,f^{-1}f_kx_k $$ separates $x$ from $P_{\mathbb{H}^2}$ (as $F$ acts by isometries,
it sends the convex polyhedron with vertices  $$f_ix_i,
f_jx_j,f_kx_k, f x, p_{\mathbb{H}^2}(f_i
x_i),p_{\mathbb{H}^2}(f_j x_j), p_{\mathbb{H}^2}(f_i
x_i), p_{\mathbb{H}^2}(f x)$$ on the convex polyhedron of vertices $$f^{-1}f_ix_i,
f^{-1}f_jx_j,f^{-1}f_kx_k,x,
p_{\mathbb{H}^2}(f^{-1}f_i
x_i),p_{\mathbb{H}^2}(f^{-1}f_j x_j),
p_{\mathbb{H}^2}(f^{-1}f_i x_i), p_{\mathbb{H}^2}(x)).$$

Then, if $x$ is projected in the triangle $(p_{\mathbb{H}^2}(f_ix_i),p_{\mathbb{H}^2}(f_j x_j),p_{\mathbb{H}^2}(f_k x_k))$, we have a condition which can be written 
\begin{equation}
\det ( f_i x_i - f_j x_j, f_i x_i-f_k x_k,\nonumber
f_i x_i-  x ) > 0.
\end{equation}
(we choose an orientation such that this $\det$ is $>0$, but 
we could take the other orientation, it will change nothing, the important fact is that the condition is open).

For each equation, the set of solution is an open half-space delimited by an affine hyperplane: the possible set of heights for the vertices  $(x_1,\ldots,x_n)$ is given by the intersection  (of an infinite number) of open half-spaces.

This is a contractible open set.  We want the vertices to stay
in the unit ball, then we must add for each vertex the condition that
its height must be (strictly) greater than $0$ and less than $(1-a_i^2-b_i^2)$,
where $a_i$ and  $b_i$ are the horizontal coordinates of each
$x_i$ (they are fixed by hypothesis): we intersect the
contractible open set with other open half-spaces, and the intersection remains a contractible open set. This set is
non-empty, as we built examples of  convex Fuchsian
polyhedra.
\end{proof}

An open set of $T_g(n)$ is parameterised by a (small) deformation
of a canonical polygon in the hyperbolic plane and a displacement
of the marked points inside this polygon. With a fixed height for
the vertices, a small displacement of a  convex Fuchsian
polyhedron (corresponding to a path in  $T_g(n)$), is always
convex (the convexity is a property preserved by a little
displacement of the vertices) and Fuchsian (by construction).

So we can endow $\mathcal{P}(n)$ with
the topology which makes it a fiber space based on $T_g(n)$, with
fibers homeomorphic to the open unit ball of $\mathbb{R}^n$:
\begin{proposition}\label{imm var hyp}
The space $\mathcal{P}(n)$ is a contractible manifold of
dimension $(6g-6+3n)$.
\end{proposition}
Because $T_g(n)$ is contractible manifold of dimension $(6g-6+2n)$, see e.g. \cite{nag}.
\paragraph{Triangulations.}

\begin{definition}
A \emph{(generalised) triangulation} of a compact  surface $S$
is a decomposition of $S$  by images by homeomorphisms of triangles of Euclidean space, with possible
identification of the edges or the vertices, such that the interiors of the faces  (resp. of the edges) are disjoint.
\end{definition}
This definition allows triangulations of the surface with only one
or two vertices. For example, take a canonical polygon such as
defined in the preceding section. Take a vertex of this polygon,
and join it with the other vertices of the polygon. By identifying
the edges of the polygon, we have a triangulation of the resulting
surface with only one vertex.

We want to know the number of edges $e$ for such a triangulation
with $n$ vertices of a compact surface of genus  $g$. As the Euler
characteristic is $\chi (g) = (2-2g)$ we have $f-e+2= 2-2g$, where
$f$ is the number of faces. As the faces are supposed to be
triangles,
 $\displaystyle f= \frac{2}{3} a$ and then
\begin{equation}
a=6g-6+3n.\nonumber
\end{equation}

\paragraph{Local description of the space of polyhedra.}

Take a subdivision of each  face of a convex Fuchsian polyhedron $P$ in
triangles (such that the resulting triangulation has no more
vertices than  the polyhedron, and is invariant
under the action of $\rho(\Gamma)$). For  such a
triangulation on $P$, we get a map $Ed_P$  which sends each
convex Fuchsian polyhedron in a neighbourhood of $P$ in $\mathcal{P}(n)$ to the
square  of the length  of  the edges of the triangulation in a
fundamental domain for the Fuchsian group action. As this triangulation
of $P$ provides a triangulation of the surface $S$, the map $Ed_P$ has its values in   $\mathbb{R}^{6g-g+3n}$.

 The map $Ed_P$ associates to  each Fuchsian polyhedron 
a set of  $(6g-g+3n)$ real numbers among all the $d_{\mathbb{H}^3}(\gamma x_i, \mu x_j)^2$, where $\gamma,\mu\in\rho(\Gamma)$, $i,j=1,\ldots,n$ and  $(x_1,\ldots,x_n)$ are the vertices of the polyhedron. It is in particular a  $C^1$  map (the description of the  topology of  $\mathcal{P}(n)$  says that 
for a neighbourhood of $P$ in $\mathcal{P}(n)$ the vertices belong to open sets of the hyperbolic space).

By the local inverse theorem, Theorem \ref{theoreme de rigidite} says exactly that  $Ed_P$ is a local homeomorphism.

\subsection{Set of metrics with conical singularities}

By standard methods involving Voronoi regions and Delaunay cellulations, 
it is known 
\cite{Rivintriangulation,Thurart1,TroTri}
 that
for each constant curvature metric with conical singularities on $S$ with constant sign singular curvature there exists a geodesic triangulation  such that the vertices of the triangulation are exactly the singular points. This allows us to see such a metric as a gluing of (geodesic) hyperbolic triangles.

(Actually, we don't need this result, because in the following we could consider only the metrics given by the induced metric 
on convex Fuchsian polyhedra. In this case, the geodesic triangulation of the metric is  given by a triangulation of the faces of the 
polyhedron).

We denote by
\begin{itemize}
\item  $\mathcal{M}(n)$ the space of Riemannian metrics on $S$ minus $n$ points. It is endowed with the following $C^k$ topology: two metrics
are close if their coefficients  until those of their $k$th derivative in any local chart are close (we don't care
which $k>2$);
\item $\widetilde{\mbox{Cone}}(n)\subset \mathcal{M}(n)$ the set of  hyperbolic metrics with $n$ conical singularities of
positive singular curvature on $S$, seen as Riemannian metrics after removing the singular points;
\item  $\mbox{Cone}(n)$ the quotient of $\widetilde{\mbox{Cone}}(n)$ by the isotopies of $S$ minus
the $n$ marked points;
\item $\widetilde{M}^T$ - where $T$ is a geodesic triangulation of an element of  $\widetilde{\mbox{Cone}}(n)$ - the set of the metrics of $\widetilde{\mbox{Cone}}(n)$ which admit a  geodesic triangulation homotopic to $T$;
\item  $\mbox{Conf}(n)$ the set of conformal structures  on $S$ with $n$
marked points.
\end{itemize}

\paragraph{Topology of the set of metrics.}

A Theorem of Mc Owen--Troyanov \cite{owen}\cite[Theorem A]{Troyanovarticle2} says that
there is a bijection between $\widetilde{\mbox{Cone}}(n)$ and $\mbox{Conf}(n)\times ]0,2\pi [^n$. More precisely, for  a conformal structure on $S$, $n$ points on $S$ and $n$ real numbers $\alpha_i$ between $0$ and $2\pi$,  there exists a unique conformal hyperbolic metric with cone singularities with angles $\alpha_i$ on $S$ (the singular curvature is $(2\pi - \alpha_i)$).

As the Teichm\"uller space $T_g(n)$ is the quotient of  $\mbox{Conf}(n)$  by the isotopies of $S$ minus its marked points,
 $\mbox{Cone}(n)$ is in bijection with $T_g(n)\times ]0,2\pi [^n$, and we endow $\mbox{Cone}(n)$ with the topology which makes this bijection a homeomorphism:
\begin{proposition}
The set $\mbox{Cone}(n)$ is a contractible manifold of dimension $(6g-6+3n)$.
\end{proposition}

\paragraph{Local description of the set of metrics.}

We call $\widetilde{Ed}_T$ the
map  from $\widetilde{M}^T$ to $\mathbb{R}^{6g-6+3n}$ which associates to each element of
$\widetilde{M}^T$ the square of the length of the edges of the triangulation. The (square of) the distance between two points of $S$ is a continuous
function from  $\mathcal{M}(n)$ to $\mathbb{R}$. 

Remark that $\widetilde{M}^T$ is non empty if we consider a metric 
given by the induced metric on a convex Fuchsian polyhedron, on which we fix a triangulation. Moreover, around a point of $\widetilde{M}^T$, $\widetilde{Ed}_T$ has its values in an open set of $\mathbb{R}^{6g-6+3n}$: if we change a little the length of the $(6g-6+3n)$ edges, the resulting metric will be still in $\widetilde{M}^T$, because the conditions to remain a hyperbolic triangle and   that the sum of the angles around each vertex don't exceed $2\pi$ are open conditions.

   Let  $i_T$ be the canonical inclusion of $\widetilde{M}^T$ (endowed with the topology of $\mathcal{M}(n)$) in $\widetilde{\mbox{Cone}}(n)$ (endowed with the topology of $\mbox{Conf}(n)\times ]0,2\pi [^n$).
    
    The composition of $i_T$ with the projection onto $\mbox{Conf}(n)$ is the map which associates to each metric its conformal structure, and it is  continuous mapping as by definition $\mbox{Conf}(n)$ is the quotient of
$\mathcal{M}(n)$ by the set of positive real-valued functions on
$S$ minus its marked points; and the composition of $i_T$ with the projection on  $]0,2\pi [^n$ is obviously continuous.

Then $i_T$ is continuous and injective: it is a local homeomorphism and then $\widetilde{Ed}_T$ is a 
continuous map on $\widetilde{M}^T\subset \widetilde{\mbox{Cone}}(n)$. Moreover, modulo the isotopies of the surface, the map $\widetilde{Ed}_T$ becomes an injective mapping
$Ed_T$ from $M^T\subset \mbox{Cone}(n)$ (the quotient of
$\widetilde{M}^T$ by the isotopies) to $\mathbb{R}^{6g-6+3n}$. And as it has its values in an open set of $\mathbb{R}^{6g-6+3n}$ and as
the dimension of $\mbox{Cone}(n)$ is $(6g-6+3n)$, $Ed_T$ is a local homeomorphism from $\mbox{Cone}(n)$ to $\mathbb{R}^{6g-6+3n}$ (defined around the images of Fuchsian polyhedra).
 
\paragraph{Realisation of metrics.}

We denote by $\mathcal{I}(n)$ the map from $\mathcal{P}(n)$ to  $\mbox{Cone}(n)$, which associates to each convex
Fuchsian polyhedral embedding in hyperbolic space its induced metric.
This map is well defined, because the induced metric on the quotient of a convex Fuchsian polyhedron (by $\rho(\Gamma)$) is isometric
to a hyperbolic metric with conical singularities with positive singular curvature on $S$. 

Let $P$ be a  convex Fuchsian polyhedron, and  $m$ its induced metric.
We consider a triangulation of $m$ given by a subdivision of the faces of $P$ in triangles.
Obviously, the (square of) the lengths of the edges of the triangulation of $P$ are the same that
the (square of) the lengths of the edges of the triangulation of $m:=\mathcal{I}(n)(P)$. It means that locally (recall that the maps $Ed_P$
and $Ed_T$ are defined only locally around $P$ and $m$):
\begin{equation} Ed_T\circ \mathcal{I}(n)\circ Ed_P^{-1}= id.\nonumber
\end{equation}
\begin{equation}
\begin{CD}
\mathcal{P}(n) @>\mathcal{I}(n)>> \mbox{Cone}(n) \\
@VEd_PVV	@VVEd_TV \\
\mathbb{R}^{6g-6+3n} @= \mathbb{R}^{6g-6+3n}
\end{CD}\nonumber
\end{equation}
 From this we deduce immediately that
$\mathcal{I}(n)$ is continuous and locally injective. Moreover, 
the map $\mathcal{I}(n)$ is proper: this will be proved in the next paragraph.

Then $\mathcal{I}(n)$
is a covering map. But as $\mathcal{P}(n)$ and  $\mbox{Cone}(n)$ are connected and simply connected,
it is a homeomorphism.

Let $\mbox{Mod}(n)$ be the quotient of the group of the homeomorphisms of $S$
minus $n$ points by its subgroup of isotopies.

Then the homeomorphism $\mathcal{I}(n)$ gives a bijection between
$\mathcal{P}(n)/\mbox{Mod}(n)$ and $\mbox{Cone}(n)/\mbox{Mod}(n)$, and this is exactly the statement of
 Theorem \ref{realisation}.

\paragraph{Properness of $\mathcal{I}(n)$.}

We will use the following characterisation of a proper map: 
$\mathcal{I}(n)$ is proper if, for each sequence $(P_k)_k$ in $\mathcal{P}(n)$ such that 
the sequence $(m_k)_k$ converges in $\mbox{Cone}(n)$ (with $m_k:=\mathcal{I}(n)(P_k)$), then 
$(P_k)_k$ converges in $\mathcal{P}(n)$ (maybe up to the extraction of a sub-sequence).

Suppose that  $(m_k)_k$ converges to $m_{\infty}\in \mbox{Cone}(n) $.  
As each $m_k$ is a hyperbolic metric with $n$ conical singularities with positive singular curvature on $S$, the convergence implies:
\begin{itemize}
\item[i)] a uniform bound on the distance between every pair of singular points of the surface for all $k$;
\item[ii)] a uniform bound on the values of the angles at the singular points for all $k$, strictly between $0$ and 
$2\pi$;
\item[iii)] a uniform bound on the lengths of the closed geodesics for all $k$;
\item[iv)] a uniform bound on the areas of the metrics for all $k$.
\end{itemize}

From this we deduce the following assertions (all supposed for a $k$ sufficiently big, and recall that we call the height of a vertex its distance from $P_{\mathbb{H}^2}$):

 \textit{If the height of one vertex of the polyhedra goes to infinity,  the heights of all the vertices go to infinity.} 
  
Suppose there exists a vertex which height doesn't go to infinity, incident to a vertex which height goes to infinity. Then 
the length of the geodesic between them goes to infinity, that is impossible by $i)$.

 \textit{The heights of all the vertices can't go to infinity.}

  By $iv)$, on each $P_k$ there is a subset $R$, homeomorphic to a closed disc,  bounding a fundamental domain for the action of $\rho_k(\Gamma)$ and which area  is fixed for all $k$. We  consider  the projection  $D_k$ of $R$ on the umbilic surface realising the minimum $d_k$ of the distance between $P_k$ and   $P_{\mathbb{H}^2}$. 
 By orthogonality of the projection, the area of $D_k$ is less than the area of $R$. The projection of $D_k$ on $P_{\mathbb{H}^2}$  is a closed set $\overline{D}_k$ bounding 
a fundamental domain for the action of  $\rho_k(\Gamma)$ on $P_{\mathbb{H}^2}$. The area of $\overline{D}_k$ is the area of $D_k$ times $(\cosh^2)^{-1} (d_k)$, then it is less than the area of $R$, which is constant, times $(\cosh^2)^{-1} (d_k)$ (recall that the induced metric on an umbilic surface at distance $d$ from $P_{\mathbb{H}^2}$ is $\cosh^2(d) \mathrm{can}_{\mathbb{H}^2}$, where $\mathrm{can}_{\mathbb{H}^2}$ is the hyperbolic metric induced on $P_{\mathbb{H}^2}$).

If the heights of all the vertices goes to infinity, it implies that
 $d_k$ goes to infinity, then the area of $\overline{D}_k$  goes to zero, and so the area of the fundamental domain for the action of  $\rho_k(\Gamma)$ in $P_{\mathbb{H}^2}$ goes to zero. This is impossible by 
 the Gauss--Bonnet Theorem: the area of a fundamental domain on  $P_{\mathbb{H}^2}$ for the action 
 of a $\rho_k(\Gamma)$ is constant for all $k$ (equal to minus the Euler characteristic of $S$).

  \textit{The lengths of the edges of the canonical polygons associated to the $P_k$ and the distances between the marked points inside it don't diverge.}  
  
  The distance between two vertices of the canonical polygon associated to  $P_k$ in 
 $P_{\mathbb{H}^2}$  is  less than the distance 
 on $P_k$ between the points which are projected onto 
 these vertices (because by orthogonality of the projection the distance in  $P_{\mathbb{H}^2}$ is smaller than the distance between the points 
 of $P_k$ in $\mathbb{H}^3$, which is itself smaller than the induced distance on $P_k$), and this distance is bounded by $i)$, so it converges (maybe after extracting a sub-sequence). 
  The same argument shows that the distance between 
  two marked points inside the canonical polygon converges. 
 
 \textit{The distance between two marked points inside a canonical polygon 
can't go to $0$.}
 
  If it occurs, it implies that the two vertices 
  $s_k$ and $s'_k$ of $P_k$ go to two points on the  geodesic 
 joining $s_k$ and $p_{\mathbb{H}^2}(s_k)$ (because we have seen that the heights of the vertices are bounded). These limit points are distinct, because the distance between two singular points is bounded. Then one is lower than the other, say $s_k$ for example. As there is an infinite number of vertices on each polyhedron  accumulating on 
  the entire boundary at infinity of $P_{\mathbb{H}^2}$, there exists at least three vertices lower than $s_k$ and such that 
 their projection on  $P_{\mathbb{H}^2}$ forms a triangle containing the projection of $s_k$: $s_k$  is contained in the interior of the convex hull of these  vertices. Then there exists a $k'$ such that for each $k>k'$, $s_k$ is in the convex hull of the others vertices: $P_k$ would be not convex, that is false. 
 
   \textit{No length of the edges of the canonical polygon goes to $0$.}
 
  It is  as above: if it occurs, there are two vertices of the canonical polygon which collapse. The corresponding points on the polyhedra 
  can't collapse, because the geodesic (on the polyhedron) between them corresponds to a closed geodesic for the induced metric on $S$ and the lengths of closed geodesics are bounded by $iii)$. Then these two points go to two points on the same 
  horizontal geodesic, and again it is in contradiction with the convexity
  of the polyhedra.\\
 
  These assertions prove that 
  the sequence of canonical polygons associated to the $P_k$ converges to 
  a compact polygon of  $P_{\mathbb{H}^2}$, with exactly $4g$ edges and $n$ marked points. As the definition of a canonical polygon shows that the limit of a converging sequence of canonical polygons is a canonical polygon, then 
   $\mathcal{S}(P_k)$ converges in $T_g(n)$. We denote by $[h]$ its limit.
   
   For $k$ sufficiently big, we can consider that all the $\mathcal{S}(P_k)$ are in a neighbourhood of $[h]$ sufficiently small to trivialise the fibration 
 $\mathcal{P}(n)$. That
 means that we can write the $P_k$ as couples $([h]_k,H_k)$, where 
 $H_k$ is the heights of the vertices of $P_k$, and we have already seen that these heights converge.
 
 Then $(P_k)_k$ converges to a Fuchsian polyhedron. It remains to check that it is convex with $n$ vertices: by convexity of the $P_k$, no  vertex converges to a point inside the convex hull of the other vertices, and no one converges to a point on the boundary of the convex hull of the other vertices, because 
the angles at the singularities of $m_{\infty}$ are all less than $2\pi$ by $ii)$.

\bibliographystyle{alpha}
\bibliography{biblio}
\end{document}

%% file: angleconvexe.pstex_t
\begin{picture}(0,0)%
\includegraphics{angleconvexe.pstex}%
\end{picture}%
\setlength{\unitlength}{4144sp}%
\begingroup\makeatletter\ifx\SetFigFont\undefined%
\gdef\SetFigFont#1#2#3#4#5{%
  \reset@font\fontsize{#1}{#2pt}%
  \fontfamily{#3}\fontseries{#4}\fontshape{#5}%
  \selectfont}%
\fi\endgroup%
\begin{picture}(4295,1876)(-7,-1025)
\put(676, 29){\makebox(0,0)[lb]{\smash{{\SetFigFont{10}{12.0}{\rmdefault}{\mddefault}{\updefault}{\color[rgb]{0,0,0}$p$}%
}}}}
\put(631,-961){\makebox(0,0)[lb]{\smash{{\SetFigFont{10}{12.0}{\rmdefault}{\mddefault}{\updefault}{\color[rgb]{0,0,0}$p_-$}%
}}}}
\put(2386,-961){\makebox(0,0)[lb]{\smash{{\SetFigFont{10}{12.0}{\rmdefault}{\mddefault}{\updefault}{\color[rgb]{0,0,0}$p_-$}%
}}}}
\put(2791,704){\makebox(0,0)[lb]{\smash{{\SetFigFont{10}{12.0}{\rmdefault}{\mddefault}{\updefault}{\color[rgb]{0,0,0}$z_{i}$}%
}}}}
\put(4141,479){\makebox(0,0)[lb]{\smash{{\SetFigFont{10}{12.0}{\rmdefault}{\mddefault}{\updefault}{\color[rgb]{0,0,0}$z_{i+1}$}%
}}}}
\put(1711,569){\makebox(0,0)[lb]{\smash{{\SetFigFont{10}{12.0}{\rmdefault}{\mddefault}{\updefault}{\color[rgb]{0,0,0}$z_{i-1}$}%
}}}}
\put(3781,-61){\makebox(0,0)[lb]{\smash{{\SetFigFont{10}{12.0}{\rmdefault}{\mddefault}{\updefault}{\color[rgb]{0,0,0}$r_{i+1}$}%
}}}}
\put(2431,-466){\makebox(0,0)[lb]{\smash{{\SetFigFont{10}{12.0}{\rmdefault}{\mddefault}{\updefault}{\color[rgb]{0,0,0}$\beta_i$}%
}}}}
\put(2746,-466){\makebox(0,0)[lb]{\smash{{\SetFigFont{10}{12.0}{\rmdefault}{\mddefault}{\updefault}{\color[rgb]{0,0,0}$\beta_{i+1}$}%
}}}}
\put(1891,-151){\makebox(0,0)[lb]{\smash{{\SetFigFont{10}{12.0}{\rmdefault}{\mddefault}{\updefault}{\color[rgb]{0,0,0}$r_{i-1}$}%
}}}}
\put(2611,209){\makebox(0,0)[lb]{\smash{{\SetFigFont{10}{12.0}{\rmdefault}{\mddefault}{\updefault}{\color[rgb]{0,0,0}$r_{i}$}%
}}}}
\put(2791,-241){\makebox(0,0)[lb]{\smash{{\SetFigFont{10}{12.0}{\rmdefault}{\mddefault}{\updefault}{\color[rgb]{0,0,0}$\alpha_i$}%
}}}}
\put(2116,614){\makebox(0,0)[lb]{\smash{{\SetFigFont{12}{14.4}{\rmdefault}{\mddefault}{\updefault}{\color[rgb]{0,0,0}$l_{i-1}$}%
}}}}
\put(3331,614){\makebox(0,0)[lb]{\smash{{\SetFigFont{12}{14.4}{\rmdefault}{\mddefault}{\updefault}{\color[rgb]{0,0,0}$l_{i}$}%
}}}}
\end{picture}%

%% file: hyperboloide.pstex_t
\begin{picture}(0,0)%
\includegraphics{hyperboloide.pstex}%
\end{picture}%
\setlength{\unitlength}{4144sp}%
\begingroup\makeatletter\ifx\SetFigFont\undefined%
\gdef\SetFigFont#1#2#3#4#5{%
  \reset@font\fontsize{#1}{#2pt}%
  \fontfamily{#3}\fontseries{#4}\fontshape{#5}%
  \selectfont}%
\fi\endgroup%
\begin{picture}(5600,3153)(-367,-2307)
\put(2235,172){\rotatebox{3.0}{\makebox(0,0)[lb]{\smash{{\SetFigFont{6}{7.2}{\rmdefault}{\mddefault}{\updefault}{\color[rgb]{0,0,0}$x$}%
}}}}}
\put(2931,245){\rotatebox{3.0}{\makebox(0,0)[lb]{\smash{{\SetFigFont{6}{7.2}{\rmdefault}{\mddefault}{\updefault}{\color[rgb]{0,0,0}$Z_l$}%
}}}}}
\put(264,433){\rotatebox{3.0}{\makebox(0,0)[lb]{\smash{{\SetFigFont{7}{8.4}{\rmdefault}{\mddefault}{\updefault}{\color[rgb]{0,0,0}$\mb{H}^2$}%
}}}}}
\put(1441,-2266){\rotatebox{3.0}{\makebox(0,0)[lb]{\smash{{\SetFigFont{6}{7.2}{\rmdefault}{\mddefault}{\updefault}{\color[rgb]{0,0,0}$0$}%
}}}}}
\put(3736,-2131){\rotatebox{3.0}{\makebox(0,0)[lb]{\smash{{\SetFigFont{6}{7.2}{\rmdefault}{\mddefault}{\updefault}{\color[rgb]{0,0,0}$0$}%
}}}}}
\put(1306,-871){\makebox(0,0)[lb]{\smash{{\SetFigFont{6}{7.2}{\rmdefault}{\mddefault}{\updefault}{\color[rgb]{0,0,0}$x_c$}%
}}}}
\put(3736,-826){\rotatebox{3.0}{\makebox(0,0)[lb]{\smash{{\SetFigFont{6}{7.2}{\rmdefault}{\mddefault}{\updefault}{\color[rgb]{0,0,0}$1$}%
}}}}}
\put(3511,524){\rotatebox{3.0}{\makebox(0,0)[lb]{\smash{{\SetFigFont{6}{7.2}{\rmdefault}{\mddefault}{\updefault}{\color[rgb]{0,0,0}$\cosh(\mu)$}%
}}}}}
\put(1036,434){\rotatebox{3.0}{\makebox(0,0)[lb]{\smash{{\SetFigFont{6}{7.2}{\rmdefault}{\mddefault}{\updefault}{\color[rgb]{0,0,0}$\cosh(\mu)$}%
}}}}}
\put(4366,614){\rotatebox{3.0}{\makebox(0,0)[lb]{\smash{{\SetFigFont{6}{7.2}{\rmdefault}{\mddefault}{\updefault}{\color[rgb]{0,0,0}$Z_l$}%
}}}}}
\put(2308,-886){\rotatebox{3.0}{\makebox(0,0)[lb]{\smash{{\SetFigFont{6}{7.2}{\rmdefault}{\mddefault}{\updefault}{\color[rgb]{0,0,0}$\Phi(Z)_l$}%
}}}}}
\put(4006,-736){\rotatebox{3.0}{\makebox(0,0)[lb]{\smash{{\SetFigFont{6}{7.2}{\rmdefault}{\mddefault}{\updefault}{\color[rgb]{0,0,0}$\Phi(Z)_l$}%
}}}}}
\end{picture}%

%% file: artrealhyp2.bbl
\begin{thebibliography}{ILTC01}

\bibitem[Ale05]{Aleks}
A.~D. Alexandrov.
\newblock {\em Convex polyhedra}.
\newblock Springer Monographs in Mathematics. Springer-Verlag, Berlin, 2005.
\newblock Translated from the 1950 Russian edition by N. S. Dairbekov, S. S.
  Kutateladze and A. B. Sossinsky, With comments and bibliography by V. A.
  Zalgaller and appendices by L. A. Shor and Yu. A. Volkov.

\bibitem[Ber77]{berger5}
Marcel Berger.
\newblock {\em G\'eom\'etrie. {V}ol. 5}.
\newblock CEDIC, Paris, 1977.
\newblock La sph\`ere pour elle-m\^eme, g\'eom\'etrie hyperbolique, l'espace
  des sph\`eres. [The sphere itself, hyperbolic geometry, the space of
  spheres].

\bibitem[Bus58]{Buse}
Herbert Busemann.
\newblock {\em Convex surfaces}.
\newblock Interscience Tracts in Pure and Applied Mathematics, no. 6.
  Interscience Publishers, Inc., New York, 1958.

\bibitem[Bus92]{Buser}
Peter Buser.
\newblock {\em Geometry and spectra of compact {R}iemann surfaces}, volume 106
  of {\em Progress in Mathematics}.
\newblock Birkh\"auser Boston Inc., Boston, MA, 1992.

\bibitem[Fil]{artrealisationlorentz}
Fran\c{c}ois Fillastre.
\newblock Polyhedral realisation of metrics with conical singularities on
  compact surfaces in {L}orentzian space-forms.
\newblock In preparation.

\bibitem[GHL90]{GaHuLa}
Sylvestre Gallot, Dominique Hulin, and Jacques Lafontaine.
\newblock {\em Riemannian geometry}.
\newblock Universitext. Springer-Verlag, Berlin, second edition, 1990.

\bibitem[Gro86]{gromovdiffrel}
Mikhael Gromov.
\newblock {\em Partial differential relations}, volume~9 of {\em Ergebnisse der
  Mathematik und ihrer Grenzgebiete (3) [Results in Mathematics and Related
  Areas (3)]}.
\newblock Springer-Verlag, Berlin, 1986.

\bibitem[ILTC01]{TroTri}
C.~Indermitte, Th.~M. Liebling, M.~Troyanov, and H.~Cl{\'e}men{\c{c}}on.
\newblock Voronoi diagrams on piecewise flat surfaces and an application to
  biological growth.
\newblock {\em Theoret. Comput. Sci.}, 263(1-2):263--274, 2001.
\newblock Combinatorics and computer science (Palaiseau, 1997).

\bibitem[Lab92]{lab2}
Fran{\c{c}}ois Labourie.
\newblock M\'etriques prescrites sur le bord des vari\'et\'es hyperboliques de
  dimension {$3$}.
\newblock {\em J. Differential Geom.}, 35(3):609--626, 1992.

\bibitem[LS00]{SchLab}
Fran{\c{c}}ois Labourie and Jean-Marc Schlenker.
\newblock Surfaces convexes fuchsiennes dans les espaces lorentziens \`a
  courbure constante.
\newblock {\em Math. Ann.}, 316(3):465--483, 2000.

\bibitem[McO88]{owen}
Robert~C. McOwen.
\newblock Point singularities and conformal metrics on {R}iemann surfaces.
\newblock {\em Proc. Amer. Math. Soc.}, 103(1):222--224, 1988.

\bibitem[Nag88]{nag}
Subhashis Nag.
\newblock {\em The complex analytic theory of {T}eichm\"uller spaces}.
\newblock Canadian Mathematical Society Series of Monographs and Advanced
  Texts. John Wiley \& Sons Inc., New York, 1988.
\newblock A Wiley-Interscience Publication.

\bibitem[O'N83]{oneill}
Barrett O'Neill.
\newblock {\em Semi-{R}iemannian geometry}, volume 103 of {\em Pure and Applied
  Mathematics}.
\newblock Academic Press Inc. [Harcourt Brace Jovanovich Publishers], New York,
  1983.
\newblock With applications to relativity.

\bibitem[Pog73]{Pogo}
A.~V. Pogorelov.
\newblock {\em Extrinsic geometry of convex surfaces}.
\newblock American Mathematical Society, Providence, R.I., 1973.
\newblock Translated from the Russian by Israel Program for Scientific
  Translations, Translations of Mathematical Monographs, Vol. 35.

\bibitem[RH93]{RivHod}
Igor Rivin and Craig~D. Hodgson.
\newblock A characterization of compact convex polyhedra in hyperbolic
  {$3$}-space.
\newblock {\em Invent. Math.}, 111(1):77--111, 1993.

\bibitem[Riv]{Rivintriangulation}
Igor Rivin.
\newblock {Extra-large metrics}.

\bibitem[Riv86]{rivinthese}
Igor Rivin.
\newblock {\em On geometry of convex polyhedra in hyperbolic 3-space}.
\newblock PhD thesis, Princeton University, June 1986.

\bibitem[Rou04]{rousset}
Mathias Rousset.
\newblock Sur la rigidit\'e de poly\`edres hyperboliques en dimension 3: cas de
  volume fini, cas hyperid\'eal cas fuchsien.
\newblock {\em Bull. Soc. Math. France}, 132(2):233--261, 2004.

\bibitem[Sab04]{sabitov}
I.~Kh. Sabitov.
\newblock Around the proof of the {L}egendre-{C}auchy lemma on convex polygons.
\newblock {\em Sibirsk. Mat. Zh.}, 45(4):892--919, 2004.

\bibitem[Sch]{Schpoly}
Jean-Marc Schlenker.
\newblock {Hyperbolic manifolds with polyhedral boundary}.
\newblock arXiv:math.GT/0111136.

\bibitem[Sch06]{Schconvex}
Jean-Marc Schlenker.
\newblock Hyperbolic manifolds with convex boundary.
\newblock {\em Invent. Math.}, 163(1):109--169, 2006.

\bibitem[Spi79]{SPI5}
Michael Spivak.
\newblock {\em A comprehensive introduction to differential geometry. {V}ol.
  {V}}.
\newblock Publish or Perish Inc., Wilmington, Del., second edition, 1979.

\bibitem[Thu98]{Thurart1}
William~P. Thurston.
\newblock Shapes of polyhedra and triangulations of the sphere.
\newblock In {\em The Epstein birthday schrift}, volume~1 of {\em Geom. Topol.
  Monogr.}, pages 511--549 (electronic). Geom. Topol. Publ., Coventry, 1998.

\bibitem[Tro91]{Troyanovarticle2}
Marc Troyanov.
\newblock Prescribing curvature on compact surfaces with conical singularities.
\newblock {\em Trans. Amer. Math. Soc.}, 324(2):793--821, 1991.

\bibitem[ZVC80]{ZVC}
Heiner Zieschang, Elmar Vogt, and Hans-Dieter Coldewey.
\newblock {\em Surfaces and planar discontinuous groups}, volume 835 of {\em
  Lecture Notes in Mathematics}.
\newblock Springer, Berlin, 1980.
\newblock Translated from the German by John Stillwell.

\end{thebibliography}
